\documentclass[final,onefignum,onetabnum]{siamonline190516}
\usepackage[scale=0.8]{geometry} 

\usepackage{lipsum}
\usepackage{amsfonts}
\usepackage{graphicx}
\usepackage{epstopdf}
\usepackage{algorithmic}
\ifpdf
  \DeclareGraphicsExtensions{.eps,.pdf,.png,.jpg}
\else
  \DeclareGraphicsExtensions{.eps}
\fi

\usepackage{datetime}
\newdateformat{monthyeardate}{%
	\monthname[\THEMONTH] \THEDAY, \THEYEAR}

\usepackage{academicons}
\usepackage{xcolor}
\renewcommand{\orcid}[1]{\href{https://orcid.org/#1}{\textcolor[HTML]{A6CE39}{orcid.org/#1}}}

\usepackage{amsmath}
\allowdisplaybreaks
\usepackage{amssymb}
\usepackage{commath}
\usepackage{mathtools}
\usepackage{bbm}
\usepackage{bm}

\usepackage{color}
\usepackage{graphicx}
\usepackage[small]{caption}
\usepackage{subcaption}

\usepackage{relsize}
\usepackage{adjustbox}
\usepackage{algorithmic}
\usepackage{booktabs}
\usepackage{tikz}
\usepackage{pifont}
\usetikzlibrary{bayesnet} 

\usepackage{enumitem}
\setlist[enumerate]{leftmargin=.5in}
\setlist[itemize]{leftmargin=.5in}

\newsiamremark{remark}{Remark}
\newsiamremark{example}{Example}
\newsiamremark{hypothesis}{Hypothesis}
\crefname{hypothesis}{Hypothesis}{Hypotheses}
\newsiamthm{claim}{Claim}

\headers{Generalized sparse Bayesian learning}{J.\ Glaubitz, A.\ Gelb, and G.\ Song}

\title{Generalized sparse Bayesian learning and application to image reconstruction
\thanks{\monthyeardate\today 
\corresponding{Jan Glaubitz} 
}}

\author{Jan Glaubitz\thanks{Department of Mathematics, Dartmouth College, Hanover, NH 03755, USA (\email{Jan.Glaubitz@Dartmouth.edu}, \orcid{0000-0002-3434-5563})}
\and 
Anne Gelb\thanks{Department of Mathematics, Dartmouth College, Hanover, NH 03755, USA (\email{Anne.E.Gelb@Dartmouth.edu}, \orcid{0000-0002-9219-4572})}
\and 
Guohui Song\thanks{Department of Mathematics and Statistics, Old Dominion University, Norfolk, VA 23529, USA (\email{gsong@odu.edu}, \orcid{0000-0002-6811-9089})}
}

\usepackage{amsopn}
\DeclareMathOperator{\diag}{diag}

\newcommand{\R}{\mathbb{R}}  
\newcommand{\C}{\mathbb{C}}
 
\newcommand{\intd}{\, \mathrm{d}} 
\renewcommand{\vec}{\, \mathrm{vec}}
\DeclareMathOperator*{\argmin}{arg\,min} 
\DeclareMathOperator*{\kernel}{kernel} 

\ifpdf
\hypersetup{
  pdftitle={Generalized SBL},
  pdfauthor={J. Glaubitz, A. Gelb, and G. Song}
}
\fi

\begin{document}

\maketitle

\begin{abstract}

Image reconstruction based on indirect, noisy, or incomplete data remains an important yet challenging task.
While methods such as compressive sensing have demonstrated high-resolution image recovery in various settings, there remain issues of robustness due to parameter tuning. 
Moreover, since the recovery is limited to a point estimate, it is impossible to quantify the uncertainty, which is often desirable.
Due to these inherent limitations, a sparse Bayesian learning approach is sometimes adopted to recover a posterior distribution of the unknown. 
Sparse Bayesian learning assumes that some linear transformation of the unknown is sparse. 
However, most of the methods developed are tailored to specific problems, with particular forward models and priors.
Here, we present a generalized approach to sparse Bayesian learning. 
It has the advantage that it can be used for various types of data acquisitions and prior information. 
Some preliminary results on image reconstruction/recovery indicate its potential use for denoising, deblurring, and magnetic resonance imaging. 
\end{abstract}




\section{Introduction} 
\label{sec:introduction} 

Many applications seek to solve the \emph{linear inverse problem} 
\begin{equation}\label{eq:intro_IP}
	\mathbf{y} = F \mathbf{x} + \boldsymbol{\nu}, 
\end{equation}
where $\mathbf{y} \in \R^m$ is a vector of indirect measurements, $\mathbf{x} \in \R^n$ is the vector of unknowns, ${F \in \R^{m \times n}}$ is a known linear forward operator, 
and $\boldsymbol{\nu} \in \R^m$ corresponds to a typically unknown noise vector  
(see \cite{groetsch1993inverse,vogel2002computational,hansen2010discrete} and references therein). 
In particular, \cref{eq:intro_IP} can be associated with signal or image reconstruction \cite{natterer2001mathematical,hansen2006deblurring,stark2013image}.
In this regard it is often reasonable to assume that some linear transformation of the unknown solution $\mathbf{x}$, say $R \mathbf{x}$, is sparse. 
A common approach is to consider the \emph{$\ell^1$-regularized inverse problem} 
\begin{equation}\label{eq:l1_RIP}
	\mathbf{x}_{\lambda} 
		= \argmin_{\mathbf{x}} \left\{ \| F \mathbf{x} - \mathbf{y} \|_2^2 + \lambda \| R \mathbf{x} \|_1 \right\},
\end{equation}
where $R \in \R^{k \times n}$ is referred to as the \emph{regularization operator} and $\lambda > 0$ as the \emph{regularization parameter}. 
The motivation for this approach is that the $\ell^1$-norm, $\|\cdot||_1$, serves as a convex surrogate for the $\ell^0$-``norm", $\|\cdot\|_0$. 
Thus, \cref{eq:l1_RIP} balances data fidelity, noise, and the sparsity assumption on $R \mathbf{x}$, while still enabling efficient computations \cite{donoho2006compressed,eldar2012compressed,foucart2017mathematical}.  
However, an often encountered difficulty for the $\ell^1$-regularized inverse problem \cref{eq:l1_RIP} is the selection of an appropriate regularization parameter $\lambda$. 
This parameter can critically influence the quality of the regularized reconstruction $\mathbf{x}_{\lambda}$ \cite{golub1979generalized,colton1997simple,hansen1999curve,sanders2020effective,lanza2020residual}. 
In part for this reason, many statistical approaches have been proposed for regularized inverse problems \cite{kaipio2006statistical,calvetti2007introduction,stark2013image}. 
Another advantage in using statistical approaches is that they may allow for uncertainty quantification in the reconstructed solution \cite{calvetti2008hypermodels}. 
For example, the hierarchical Bayesian formulation of the $\ell^1$-regularized inverse problem \cref{eq:l1_RIP}, see \cite{kaipio2006statistical,calvetti2007introduction}, is based on extending $\mathbf{x}$, $\mathbf{y}$, and all other involved parameters, which we collectively write as $\boldsymbol{\theta}$, into random variables. 
Consequently $\mathbf{x}$, $\mathbf{y}$, and $\boldsymbol{\theta}$ are characterized by certain density functions, as are their relationships to each other. 
In particular, one usually considers the following density functions: 
\begin{itemize} 
    \item 
    The \emph{likelihood} $p(\mathbf{y} | \mathbf{x}, \boldsymbol{\theta})$, which is the probability density function for $\mathbf{y}$ given $\mathbf{x}$ and $\boldsymbol{\theta}$.

	\item 
    The \emph{prior} $p(\mathbf{x} | \boldsymbol{\theta})$, which is the density function for $\mathbf{x}$ given $\boldsymbol{\theta}$. 
	
	\item 
    The \emph{hyper-prior} $p(\boldsymbol{\theta})$, which is the probability density function for the parameters $\boldsymbol{\theta}$. 
	
	\item 
	The \emph{posterior} $p(\mathbf{x}, \boldsymbol{\theta} | \mathbf{y})$, which is the probability density function for the solution $\mathbf{x}$ and the parameters $\boldsymbol{\theta}$ given the data $\mathbf{y}$.
	
\end{itemize}
One can use Bayes' theorem to obtain
\begin{equation}\label{eq:Bayes_formula} 
	p(\mathbf{x}, \boldsymbol{\theta} | \mathbf{y}) 
		\propto p(\mathbf{y} | \mathbf{x}, \boldsymbol{\theta}) p(\mathbf{x}|\boldsymbol{\theta}) p(\boldsymbol{\theta}),
\end{equation} 
where ``$\propto$" means that the two sides of \cref{eq:Bayes_formula} are equal to each other up to a multiplicative constant that does not depend on $\mathbf{x}$ or $\boldsymbol{\theta}$. 
Note that the parameters $\boldsymbol{\theta}$ are now part of the problem and are no longer selected a priori. 
Furthermore, using an appropriate method for Bayesian inference allows to quantify uncertainty in the reconstructed solution $\mathbf{x}$.

There are a variety of sparsity-promoting priors to choose from, including but not limited to Laplace priors \cite{figueiredo2007majorization}, TV-priors \cite{kaipio2000statistical,babacan2008parameter}, mixture-of-Gaussian priors \cite{fergus2006removing}, and hyper-Laplacian distributions based on $\ell^p$-quasinorms with $0<p<1$ \cite{levin2007image,krishnan2009fast}. 
In this investigation we consider the well-known class of conditionally Gaussian priors given by
\begin{equation}\label{eq:Gaussian_prior} 
	p(\mathbf{x} | \boldsymbol{\beta}) 
		\propto \det(B)^{1/2} \exp\left\{ -\frac{1}{2} \mathbf{x}^T R^T B R \mathbf{x} \right\},
\end{equation}
where $B = \diag(\beta_1,\dots,\beta_k)$ is a diagonal inverse covariance matrix. 
Ideas discussed in \cite{tipping2001sparse,wipf2004sparse,chantas2006bayesian,calvetti2007gaussian,chantas2008variational,bardsley2010hierarchical,babacan2010sparse,calvetti2019hierachical,churchill2019detecting,churchill2020estimation} suggest that conditionally Gaussian priors of the form \cref{eq:Gaussian_prior} are particularly suited to promote sparsity of $R \mathbf{x}$. 
For example, the model proposed in \cite{tipping2001sparse} is designed to recover sparse representations of kernel approximations, coining the term \emph{sparse Bayesian learning (SBL)}. 
Promoting sparse solutions, as done in \cite{tipping2001sparse}, corresponds to using $R = I \in \R^{n \times n}$ as a regularization operator in \cref{eq:Gaussian_prior}. 
Further investigations that made use of SBL to promote sparse solutions include \cite{wipf2004sparse,zhang2011sparse,zhang2014spatiotemporal,calvetti2019hierachical,churchill2019detecting}.
In many applications, however, it is some linear transformation $R \mathbf{x}$ that is desired to be sparse. 
For example, total variation (TV) regularization is of particular interest in image recovery.  
Extensions of SBL for this setting have been proposed in \cite{chantas2006bayesian,calvetti2007gaussian,calvetti2008hypermodels,chantas2008variational,bardsley2010hierarchical,babacan2010sparse,churchill2020estimation}. 
That said, since the TV-regularization operator $R$ is singular, the prior  \cref{eq:Gaussian_prior} is improper. 
This prohibits the application of many of the existing SBL approaches. 
An often encountered idea therefore is to make $R \in \R^{k\times n}$ with $k < n$ invertible by introducing additional rows that are consistent with assumptions about the underlying solution. 
For example, in \cite{calvetti2007gaussian,calvetti2008hypermodels,bardsley2010hierarchical} the additional rows encode certain boundary conditions. 
The same technique can be extended to higher-order TV-regularization \cite{churchill2020estimation}. 
Unfortunately, such additional information might not always be available or may be complicated to incorporate, especially in two or more dimensions. 
Further, different types of regularization operators must be adapted on a case-by-case basis, and the resulting prior may promote undesired artificial features in the solution when the regularization operator is not carefully modified.  
The approach in \cite{chantas2006bayesian,chantas2008variational}, by contrast, depends on the assumption of  a ``commuting property" of the form $F R = R F$. 
Requiring such a commuting property is often unrealistic in  applications, however.\footnote{The dimensions of $F$ and $R$ are typically not consistent.}

\subsection*{Our Contribution}

To address these issues, we present a generalized approach to SBL for ``almost" general forward and regularization operators, $F$ and $R$. 
By ``almost'' general, we mean that the only restriction on $F$ and $R$ is that their common kernel should be trivial, $\kernel(F) \cap \kernel(R) = \{\mathbf{0}\}$, a standard assumption in regularized inverse problems \cite{kaipio2006statistical}. 
We propose an efficient numerical method for Bayesian inference that yields a full conditional posterior density $p(\mathbf{x} | \mathbf{y})$, rather than a simple point estimate, which  allows for uncertainty quantification in the solution $\mathbf{x}$. 
The present work implies that SBL can be applied to a broader class of problems than currently known.
In particular, some preliminary results on signal and image reconstruction indicate its potential use for denoising, deblurring, and magnetic resonance imaging.

\subsection*{Outline} 

The rest of this paper is organized as follows. 
\Cref{sec:Bayesian} provides some details on the sparsity promoting hierarchical Bayesian model under consideration. 
In \cref{sec:BI}, we propose an efficient numerical method for Bayesian inference. 
A series of numerical examples is presented in \cref{sec:numerics} to illustrate the descriptive span of the hierarchical Bayesian model. 
Finally, in \cref{sec:summary}, we provide some concluding thoughts. 
 

\section{The hierarchical Bayesian model} 
\label{sec:Bayesian} 

\begin{figure}[tb]
	\centering
	\resizebox{0.6\textwidth}{!}{%
	\begin{tikzpicture}
    		\node[obs] (y) {$\mathbf{y}$} ; %
		\node[latent, below=of y] (alpha) {$\boldsymbol{\alpha}$} ; %
		\node[const, left=1.0 of alpha, yshift=0.5cm] (c_noise) {$c$ \ } ; %
  		\node[const, left=1.0 of alpha, yshift=-0.5cm]  (d_noise) {$d$ \ } ; %
		\node[latent, right=1.25 of y] (x) {$\mathbf{x}$} ; %
		\node[latent, right=1.25 of x] (beta) {$\boldsymbol{\beta}$} ; %
		\node[const, right=1.0 of beta, yshift=0.5cm] (c_prior) { \ $c$} ; %
  		\node[const, right=1.0 of beta, yshift=-0.5cm]  (d_prior) {\ $d$} ; %
		\edge {alpha} {y} ; %
		\edge {x} {y} ; %
		\edge {c_noise} {alpha}; 
		\edge {d_noise} {alpha}; 
		\edge {beta} {x}; 
		\edge {c_prior} {beta}; 
		\edge {d_prior} {beta}; 
		\plate[inner sep=0.4cm, xshift=0cm, yshift=0cm] {plate1} {(y) (alpha)} {$m$}; %
		\plate[inner sep=0.4cm, xshift=0cm, yshift=0cm] {plate2} {(x)} {$n$}; %
		\plate[inner sep=0.4cm, xshift=0cm, yshift=0cm] {plate3} {(beta)} {$k$}; %
  \end{tikzpicture} 
  }%
  \caption{
  	Graphical representation of the hierarchical Bayesian model. 
  	Nodes denoted within circles correspond to random variables, while nodes without a circle correspond to parameters. 
  	Shaded circles represent observed random variables, while plain circles represent hidden random variables. 
  }
  \label{fig:graphical_model}
\end{figure}
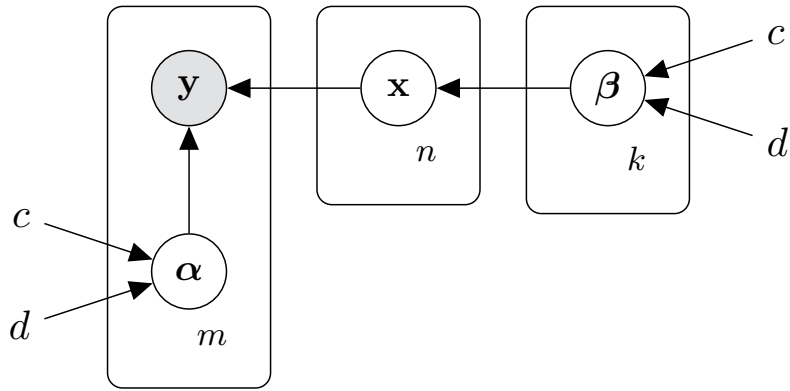

We begin by reviewing the generalized hierarchical Bayesian model considered here, which is illustrated in \cref{fig:graphical_model}. 

\subsection{The likelihood} 
\label{sub:likelihood}

The likelihood function $p(\mathbf{y} | \mathbf{x}, \boldsymbol{\alpha})$ models the connection between the solution $\mathbf{x}$, the noise parameters $\boldsymbol{\alpha}$, and the indirect measurements $\mathbf{y}$. 
It is often assumed that $\boldsymbol{\nu} \in \R^m$ in \cref{eq:intro_IP} is zero-mean i.\,i.\,d.\ normal noise with inverse variance $\alpha > 0$, that is, $\boldsymbol{\nu} \sim \mathcal{N}(\mathbf{0},\alpha^{-1} I)$.
This assumption yields the conditionally Gaussian likelihood function 
\begin{equation}\label{eq:likelihood_iid}
	p(\mathbf{y} | \mathbf{x}, \alpha) 
		= (2\pi)^{-m/2} \alpha^{m/2} \exp \left\{ -\frac{\alpha}{2} \| F\mathbf{x} - \mathbf{y} \|_2^2 \right\}.
\end{equation} 
The likelihood function given by \cref{eq:likelihood_iid} was considered, for instance, in \cite{tipping2001sparse,chantas2006bayesian,babacan2010sparse,bardsley2012mcmc,churchill2020estimation}.
By contrast, we restrict $\boldsymbol{\nu}$ to be independent but \emph{not} necessarily identically distributed. 
This translates into ${\boldsymbol{\nu} \sim \mathcal{N}(\mathbf{0},A^{-1})}$ with diagonal positive definite \emph{inverse noise covariance matrix}  
\begin{equation}\label{eq:inv_cov_matrix}
    A = \diag(\boldsymbol{\alpha}), \quad \boldsymbol{\alpha} = [\alpha_1,\dots,\alpha_m].
\end{equation}
The linear data model \cref{eq:intro_IP} then yields the generalized likelihood function 
\begin{equation}\label{eq:likelihood}
	p(\mathbf{y} | \mathbf{x}, \boldsymbol{\alpha}) 
		= (2\pi)^{-m/2} \det(A)^{1/2} \exp \left\{ - \frac{1}{2} ( F\mathbf{x} - \mathbf{y} )^T A ( F\mathbf{x} - \mathbf{y} ) \right\},
\end{equation} 
which reduces to \cref{eq:likelihood_iid} if the inverse variances $\alpha_1, \dots,\alpha_m$ are all equal to $\alpha$. 
We note that conditionally Gaussian likelihoods of the form \cref{eq:likelihood} were considered in \cite[Example 3.4]{calvetti2010subjective} in combination with smoothness promoting priors to address data outliers. 
\cref{expl:data_fusion} below motivates the weaker assumption ${\boldsymbol{\nu} \sim \mathcal{N}(\mathbf{0},A^{-1})}$ for sparsity promoting priors in the context of data fusion \cite{gustafsson2010statistical} and multi-sensor acquisition systems \cite{guerquin2011realistic,chun2017compressed}. 

\begin{example}\label{expl:data_fusion}
    Assume we have a collection of measurements $\mathbf{y}^{(d)} \in \R^{m_d}$, $d=1,\dots,D$, generated from the same source $\mathbf{x}$ from $D$ different sensors. 
    The corresponding data models are 
    \begin{equation}\label{eq:data_model_seperate}  
        \mathbf{y}^{(d)} = F^{(d)} \mathbf{x} + \boldsymbol{\nu}^{(d)}, 
        \quad d=1,\dots,D.
    \end{equation}
    Further, assume that the noise in the measurements from the same sensor is i.\,i.\,d., that is,  $\boldsymbol{\nu}^{(d)} \sim \mathcal{N}(\mathbf{0},\alpha_1^{-1}I)$. 
    However, the noise variance might differ from sensor to sensor, so that $\alpha_1 \neq \dots \neq \alpha_D$. 
    If we combine the different measurements and consider the joint data model 
    \begin{equation}\label{eq:data_model_fusion} 
        \underbrace{%
        \begin{bmatrix} \mathbf{y}^{(1)} \\ \vdots \\ \mathbf{y}^{(D)} \end{bmatrix}%
        }_{= \mathbf{y}}
        = 
        \underbrace{%
        \begin{bmatrix} F^{(1)} \\ \vdots \\ F^{(D)} \end{bmatrix}%
        }_{= F}
        \mathbf{x}
        + 
        \underbrace{%
        \begin{bmatrix} \boldsymbol{\nu}^{(1)} \\ \vdots \\ \boldsymbol{\nu}^{(D)} \end{bmatrix}%
        }_{= \boldsymbol{\nu}},
    \end{equation} 
    the stacked noise vector $\boldsymbol{\nu}$ cannot be assumed to be i.\,i.\,d., which we cannot appropriately model using the likelihood function \cref{eq:likelihood_iid}.
    However, using the more general likelihood function \cref{eq:likelihood}, we can model \cref{eq:data_model_fusion} by choosing a diagonal inverse noise covariance matrix $A$ of the form 
    \begin{equation}\label{eq:covariancematrix} 
        A = \diag(\alpha_1 I_1,\dots,\alpha_D I_D), 
    \end{equation} 
    where $I_d \in \R^{m_d \times m_d}$, $d=1,\dots,D$, denotes the identity matrix with dimensions matching the number of measurements provided by the $d$th sensor. 
\end{example}

\begin{remark}
	We note that in \cite{zhang2011clarify} it was pointed out that for classical SBL algorithms, even when the exact inverse noise variance $\alpha$ (or $A$) is known, using this fixed value instead of a variable Gamma hyper-prior can yield suboptimal reconstructions.   
\end{remark}

\subsection{The prior} 

The prior function $p(\mathbf{x} | \boldsymbol{\beta})$ models our prior belief about the unknown solution $\mathbf{x}$. 
Assume that some linear transformation of $\mathbf{x}$, say $R\mathbf{x}$, is sparse. 
The SBL approach promotes this assumed sparsity by using a conditionally Gaussian prior function, 
\begin{equation}\label{eq:prior} 
	p(\mathbf{x} | \boldsymbol{\beta}) 
		= \det(B)^{1/2} \exp\left\{ -\frac{1}{2} \mathbf{x}^T R^T B R \mathbf{x} \right\}, 
\end{equation}%
where $B = \diag(\beta_1,\dots,\beta_k)$ is referred to as the \emph{inverse prior convariance matrix}. 
See \cite{tipping2001sparse,wipf2004sparse,chantas2006bayesian,calvetti2007gaussian,chantas2008variational,bardsley2010hierarchical,babacan2010sparse,calvetti2019hierachical,churchill2019detecting,churchill2020estimation} and references therein. 
The conditionally Gaussian prior \cref{eq:prior} can be justified by its asymptotic behavior \cite{calvetti2007gaussian}. 
If we assume that the inverse variances $\beta_1,\dots,\beta_k$ are all equal, then \cref{eq:prior} favors solutions $\mathbf{x}$ for which $R \mathbf{x}$ is equal or close to zero,\footnote{For this prior, $R \mathbf{x}$ being close to zero means that $R \mathbf{x}$ has a small (unweighted) $\ell^2$-norm, $\|R \mathbf{x}\|_2$.} since these solutions have a higher probability. 
For instance, when $R\mathbf{x}$ corresponds to the total variation of $\mathbf{x}$, $[R\mathbf{x}]_j = x_{j+1} - x_j$, then \cref{eq:prior} promotes solutions $\mathbf{x}$ that have no or little variation.
However, if one of the inverse variances, say $\beta_j$, is significantly smaller than the remaining ones, a jump between $x_j$ and $x_{j+1}$ becomes more likely. 
In this way, \cref{eq:prior} promotes sparsity of $R\mathbf{x}$. 

\begin{remark}[Improper priors] 
    If $\kernel{R} \neq \{\mathbf{0}\}$, then $R^T B R$ is singular and \cref{eq:prior} becomes an improper prior.
    Most existing SBL algorithms are infeasible in this case, thus motivating us to propose an alternative method in \cref{sec:BI}. 
    In particular, the resulting difficulties for the evidence approach, which was used in the original investigation \cite{tipping2001sparse} and later in \cite{babacan2010sparse}, are addressed in \cref{app:improper_evidences}. 
\end{remark}

\subsection{The hyper-prior}

From the discussion above it is evident that the inverse variances $\beta_1,\dots,\beta_k$ must be allowed to have distinctly different values for the conditionally Gaussian prior \cref{eq:prior} to promote sparsity of $R \mathbf{x}$. 
This can be achieved by treating $\beta_1,\dots,\beta_k$ as random variables with uninformative density functions. 
A common choice is the gamma distribution with probability density function 
\begin{equation}\label{eq:pdf_gamma}
    \Gamma( x | c,d ) 
        = \frac{d^c}{\Gamma(c)} x^{c-1} e^{-dx},
\end{equation}
where $c$ and $d$ are positive shape and rate parameters. 
Furthermore, $\Gamma(\cdot)$ on the right-hand side of \cref{eq:pdf_gamma} denotes the usual gamma function \cite{artin2015gamma}. 
Note that a gamma-distributed random variable, $X \sim \Gamma(c,d)$, respectively has mean $E[X] = c/d$ and variance $V[X] = c/d^2$. 
In particular, $c \to 1$ and $d \to 0$ implies $E[X],V[X] \to \infty$, making \cref{eq:pdf_gamma} an uninformative prior. 
We therefore choose the inverse noise and prior variances, $\boldsymbol{\alpha}$ and $\boldsymbol{\beta}$, to be gamma-distributed: 
\begin{equation}\label{eq:hyper_priors} 
\begin{aligned} 
	p(\alpha_i) & = \Gamma \left( \alpha_i | c, d \right), \quad 
		i=1,\dots,m, \\ 
	p(\beta_j) & = \Gamma \left( \beta_j | c, d \right), \quad 
		j=1,\dots,k. 
\end{aligned}
\end{equation} 
By setting $c=1$ and $d \approx 0$, $\boldsymbol{\alpha}$ and $\boldsymbol{\beta}$ are free from the moderating influence of the hyper-prior and allowed to ``vary wildly" following the data. 
In our numerical tests we used $d=10^{-4}$ for all one-dimensional problems (signals) and $d=10^{-2}$ for all two-dimensional problems (images), which is similar to the choices in \cite{tipping2001sparse,babacan2010sparse,bardsley2012mcmc}. 
Future investigations will elaborate on the influence of these parameters. 
A few remarks are in order.

\begin{remark}[Conjugate hyper-priors]
   Choosing the hyper-priors $p(\alpha_i)$ and $p(\beta_j)$ as gamma distributions is convenient since the gamma distribution is a conjugate\footnote{Recall that $p(\theta)$ is a \emph{conjugate} for $p(z|\theta)$ if the posterior $p(\theta|z)$ is in the same class of densities (in this case corresponding to gamma distributions) as $p(\theta)$.} (see \cite{gelman1995bayesian,fink1997compendium,murphy2007conjugate}) for the conditionally Gaussian distributions \cref{eq:likelihood} and \cref{eq:prior}. 
\end{remark}

\begin{remark}[Informative hyper-priors]
    For simplicity we use the same hyper-prior $\Gamma(\cdot|c,d)$ and parameters $c,d$ for all components of $\boldsymbol{\alpha}$, $\boldsymbol{\beta}$. 
    If one has a reasonable a priori notion of what $\boldsymbol{\alpha}$ or $\boldsymbol{\beta}$ should be, the choice for hyper-prior could be modified correspondingly \cite{bardsley2012mcmc,calvetti2019hierachical}. 
\end{remark} 

\begin{remark}[Generalized gamma hyper-priors] 
    The use of generalized gamma distributions was recently investigated in \cite{calvetti2020sparse} and merged into a hybrid solver in \cite{calvetti2020sparsity}. 
    Although generalized gamma hyper-priors were demonstrated to promote sparsity more strongly in some cases, to not exceed the scope of the present work, we limit our discussion to usual gamma hyper-priors. 
\end{remark}
  

\section{Bayesian inference} 
\label{sec:BI} 

We now propose a Bayesian inference method for the generalized hierarchical Bayesian model from \cref{sec:Bayesian}.

\subsection{Preliminary observations}

The conditionally Gaussian prior \cref{eq:prior} and the gamma hyper-prior \cref{eq:pdf_gamma} were intentionally chosen because of their conditional conjugacy relationship. 
Some especially important implications include the following (see \cite{gelman1995bayesian}): 
\begin{align}
    p( \mathbf{y} | \mathbf{x}, \boldsymbol{\alpha} ) p(\mathbf{x} |  \boldsymbol{\beta})
	    & \propto \mathcal{N}( \mathbf{x} | \boldsymbol{\mu} ,C), \label{eq:prod_likelihood_prior} \\ 
	p( \mathbf{y} | \mathbf{x}, \boldsymbol{\alpha} ) p(\boldsymbol{\alpha}) 
	    & \propto \prod_{i=1}^m \Gamma( \alpha_i | 1/2 + c, [F\mathbf{x} - \mathbf{y}]_i^2/2 + d ), \label{eq:prod_likelihood_hyper} \\ 
    p( \mathbf{x} | \boldsymbol{\beta} ) p(\boldsymbol{\beta}) 
	    & \propto \prod_{j=1}^k \Gamma( \beta_j | 1/2 + c, [R\mathbf{x}]_j^2/2 + d ).  \label{eq:prod_prior_hyper}
\end{align}
Here the covariance matrix $C$ and the mean $\boldsymbol{\mu}$ in \cref{eq:prod_likelihood_prior} are respectively given by 
\begin{equation}\label{eq:mean_covariance2}
\begin{aligned} 
	C = \left( F^T A F + R^T B R \right)^{-1}, \quad 
	\boldsymbol{\mu} = C F^T A \mathbf{y}, 
\end{aligned}
\end{equation} 
$[F\mathbf{x} - \mathbf{y}]_i$ denotes the $i$th entry of the vector $F\mathbf{x} - \mathbf{y} \in \R^m$, and $[R\mathbf{x}]_j$ denotes the $j$th entry of the vector $R\mathbf{x} \in \R^k$.
Note that the two sides of \cref{eq:prod_likelihood_prior}, \cref{eq:prod_likelihood_hyper}, and \cref{eq:prod_prior_hyper} are equal up to a multiplicative constant that does not depend on $\mathbf{x}$, $\boldsymbol{\alpha}$, and $\boldsymbol{\beta}$, respectively. 
Finally, we stress that \cref{eq:prod_likelihood_prior} only holds if the forward operator $F \in \R^{m \times n}$ and the regularization operator $R \in \R^{k \times n}$ satisfy the \emph{common kernel condition}: 
\begin{equation}\label{eq:common_kernel}
    \kernel(F) \cap \kernel(R) = \{ \mathbf{0} \},
\end{equation} 
which is a standard assumption in regularized inverse problems \cite{kaipio2006statistical,tikhonov2013numerical}. 
Indeed, \cref{eq:common_kernel} can be interpreted as the prior (regularization) introducing a sufficient amount of complementary information to the likelihood (the given measurements) to make the problem well-posed. 
This indicates that the hierarchical Bayesian model proposed in \cref{sec:Bayesian} does not require $R$ to be invertible as long as \cref{eq:common_kernel} is satisfied.

\subsection{Proposed method: Bayesian coordinate descent} 

We are now in a position to formulate a Bayesian inference method for the generalized hierarchical Bayesian model from \cref{sec:Bayesian}. 
This method is motivated by the coordinate descent approaches \cite{friedman2010regularization,wright2015coordinate} and solves for a descriptive quantity (mode, mean, variance, etc.)~of the posterior density function $p(\mathbf{x},\boldsymbol{\alpha},\boldsymbol{\beta}|\mathbf{y})$ by alternatingly updating this quantity for $\mathbf{x}$, $\boldsymbol{\alpha}$, and $\boldsymbol{\beta}$. 
Henceforth we refer to this method as the \emph{Bayesian coordinate descent (BCD) algorithm }. 

Assume that we are interested in the expected value (mean) of the posterior, $E[\mathbf{x},\boldsymbol{\alpha},\boldsymbol{\beta}|\mathbf{y}]$. 
The BCD algorithm for this case is described in \cref{algo:BCD_mean}.

\begin{algorithm}[h!]
\caption{BCD algorithm for the mean}\label{algo:BCD_mean} 
\begin{algorithmic}[1]
    \STATE{Initialize $\boldsymbol{\alpha}^0$, $\boldsymbol{\beta}^0$, and $l=0$ } 
    \REPEAT
    	\STATE{Update $\mathbf{x}$ by setting $\mathbf{x}^{l+1} = E[ \mathbf{x} | \boldsymbol{\alpha}^{l},\boldsymbol{\beta}^{l}, \mathbf{y} ]$}
		\STATE{Update $\boldsymbol{\alpha}$ by setting $\boldsymbol{\alpha}^{l+1} = E[ \boldsymbol{\alpha} | \mathbf{x}^{l+1}, \boldsymbol{\beta}^{l}, \mathbf{y} ]$} 
		\STATE{Update $\boldsymbol{\beta}$ by setting $\boldsymbol{\beta}^{l+1} = E[ \boldsymbol{\beta} | \mathbf{x}^{l+1}, \boldsymbol{\alpha}^{l+1}, \mathbf{y} ]$}
		\STATE{Increase $l \to l+1$}
    \UNTIL{convergence or maximum number of iterations is reached}
\end{algorithmic}
\end{algorithm} 

In \cref{algo:BCD_mean} and henceforth, all variables with superscripts are treated as fixed parameters. 
That is, the expected values in \cref{algo:BCD_mean} are respectively computed w.\,r.\,t.\ $\mathbf{x}$, $\boldsymbol{\alpha}$, and $\boldsymbol{\beta}$. 
\cref{algo:BCD_mean} is simple to implement because of the particular decomposition of the posterior density function $p(\mathbf{x},\boldsymbol{\alpha},\boldsymbol{\beta}|\mathbf{y})$ provided by Bayes' theorem (see \cref{eq:Bayes_formula}): 
\begin{equation}\label{eq:Bayes_formula_BI} 
	p(\mathbf{x}, \boldsymbol{\alpha}, \boldsymbol{\beta} | \mathbf{y}) 
		\propto p(\mathbf{y} | \mathbf{x}, \boldsymbol{\alpha}) p(\mathbf{x}|\boldsymbol{\beta}) p(\boldsymbol{\alpha}) p(\boldsymbol{\beta}) 
\end{equation} 
By \cref{eq:prod_likelihood_prior,eq:prod_likelihood_hyper,eq:prod_prior_hyper}, we therefore have 
\begin{align}
    p(\mathbf{x} | \boldsymbol{\alpha}^{l},\boldsymbol{\beta}^{l}, \mathbf{y}) 
        & \propto p( \mathbf{y} | \mathbf{x}, \boldsymbol{\alpha}^{l} ) p(\mathbf{x} |  \boldsymbol{\beta}^{l})
	    \propto \mathcal{N}( \mathbf{x} | \boldsymbol{\mu} ,C), \label{eq:post1} \\ 
	p(\boldsymbol{\alpha} | \mathbf{x}^{l+1}, \boldsymbol{\beta}^{l}, \mathbf{y}) 
        & \propto p( \mathbf{y} | \mathbf{x}^{l+1}, \boldsymbol{\alpha} ) p(\boldsymbol{\alpha}) 
	    \propto \prod_{i=1}^m \Gamma( \alpha_i | 1/2 + c, [F\mathbf{x}^{l+1} - \mathbf{y}]_i^2/2 + d ), \label{eq:post2} \\ 
    p(\boldsymbol{\beta} | \mathbf{x}^{l+1}, \boldsymbol{\alpha}^{l+1}, \mathbf{y}) 
        & \propto p( \mathbf{x}^{l+1} | \boldsymbol{\beta} ) p(\boldsymbol{\beta}) 
	    \propto \prod_{j=1}^k \Gamma( \beta_j | 1/2 + c, [R\mathbf{x}^{l+1}]_j^2/2 + d ),  \label{eq:post3}
\end{align} 
where the covariance matrix $C$ and the mean $\boldsymbol{\mu}$ in \cref{eq:post1} are given as in \cref{eq:mean_covariance2} with $A = \diag(\boldsymbol{\alpha}^{l})$ and $B = \diag(\boldsymbol{\beta}^{l})$. 
Thus, the update step for $\mathbf{x}$ in \cref{algo:BCD_mean} reduces to solving the linear system 
\begin{equation}\label{eq:update_x}
    \left( F^T A F + R^T B R \right) \mathbf{x}^{l+1} = F^T A \mathbf{y}
\end{equation} 
for the mean $\mathbf{x}^{l+1}$, and the subsequent update steps for $\boldsymbol{\alpha}$ and $\boldsymbol{\beta}$ yield  
\begin{align}
    \alpha_i^{l+1} 
        & = \frac{1+2c}{\left[ F\mathbf{x}^{l+1} - \mathbf{y} \right]_i^2 + 2d}, 
        \quad i=1,\dots,m, \label{eq:update_alpha} \\ 
    \beta_j^{l+1} 
        & = \frac{1+2c}{\left[ R\mathbf{x}^{l+1} \right]_j^2 + 2d}, 
        \quad j=1,\dots,k, \label{eq:update_beta}
\end{align} 
respectively. 
Hence, \cref{algo:BCD_mean} consists of alternating between \cref{eq:update_x,eq:update_alpha,eq:update_beta}. 

\begin{remark}
    For i.\,i.\,d.\ noise, that is, the likelihood function is \cref{eq:likelihood_iid} rather than \cref{eq:likelihood}, the linear system \cref{eq:update_x} will be simplified to 
    \begin{equation}\label{eq:update_x_iid}
        \left( \alpha F^T F + R^T B R \right) \mathbf{x}^{l+1} = \alpha F^T \mathbf{y},
    \end{equation} 
    and the update step \cref{eq:update_alpha} correspondingly reduces to 
    \begin{equation}\label{eq:update_alpha_iid}
        \alpha^{l+1} = \frac{m + 2c}{\|F\mathbf{x}^{l+1} - \mathbf{y}\|_2^2 + 2d}.
    \end{equation}
\end{remark}

\begin{remark} 
	It was demonstrated in \cite{wipf2004sparse} that the cost function of classic SBL, which can be recovered from the generalized model in \cref{sec:Bayesian} for $R = I$, is non-convex with potentially many local minima that are achieved at a sparse solution. 
	Further, the cost function has a global minimum that can produce the maximally sparse solution at the posterior mean and the classic SBL algorithm based on evidence maximization is globally convergent. 
	While we numerically observed similar properties in the context of GSBL and other regularization operators $R$ (with sparsity holding for $R \mathbf{x}$ instead of $\mathbf{x}$), a detailed analysis exceeds the scope of the present paper. 
\end{remark}

\subsection{Efficient implementation of the $\mathbf{x}$-update}
\label{sub:efficient_xupdate}

If the common kernel condition \cref{eq:common_kernel} is satisfied, then the coefficient matrix on the left-hand side of \cref{eq:update_x} is symmetric and positive definite (SPD). 
For sufficiently small problems, \cref{eq:update_x} can therefore be solved efficiently using a preconditioned conjugate gradient (PCG) method \cite{saad2003iterative}. 
However, the coefficient matrix may become prohibitively large in some cases. 
To avoid any potential storage and computational issues, we implemented our method using gradient descent for the imaging problems described in \cref{sec:numerics}. 

Let $G = F^T A F + R^T B R$ and $\mathbf{b} = F^T A \mathbf{y}$ be the SPD coefficient matrix and the right-hand side of the linear system \cref{eq:update_x}, respectively.  
The solution of \cref{eq:update_x} then corresponds to the unique minimizer of the quadratic functional 
\begin{equation}\label{eq:functional}
	J(\mathbf{x}) = \mathbf{x}^T G \mathbf{x} - 2 \mathbf{x}^T \mathbf{b} 
	\quad \text{with} \quad 
	\nabla J(\mathbf{x}) = 2 \left( G \mathbf{x} - \mathbf{b} \right).
\end{equation} 
For this functional, line search minimization can be performed analytically to find the locally optimal step size $\gamma$ in every iteration. 
This allows us to use the classical gradient descent method described in \cref{algo:grad_desc} to approximate the solution $\mathbf{x}^{l+1}$ of \cref{eq:update_x}. 

\begin{algorithm}[h!]
\caption{Gradient descent method}
\label{algo:grad_desc}
\begin{algorithmic}[1]
    \STATE{Set $\mathbf{r} = \mathbf{b} - G \mathbf{x}$} 
    \REPEAT 
    		\STATE{Compute $G \mathbf{r}$ according to \cref{eq:G_effic}}
    		\STATE{Compute the step size: $\gamma = \mathbf{r}^T \mathbf{r} / \mathbf{r}^T G \mathbf{r}$} 
		\STATE{Update the solution: $\mathbf{x} + \gamma \mathbf{r}$} 
		\STATE{Update the difference: $\mathbf{r} = \mathbf{r} - \gamma G \mathbf{r}$} 
    \UNTIL{convergence or maximum number of iterations is reached}
\end{algorithmic}
\end{algorithm} 

It is important to note that the gradient in \cref{eq:functional} can be computed efficiently and without having to store the whole coefficient matrix $G$, which might be prohibitively large. 
To show this, assume that the unknown solution $\mathbf{x} \in \R^{n^2}$ corresponds to a vectorized matrix $X \in \R^{n \times n}$ and that the forward operator $F$ corresponds to applying the same one-dimensional forward operator $F_1$ to the matrix $X$ in $x$- and $y$-direction: 
\begin{equation}
    F \mathbf{x} = \mathbf{y} \iff F_1 X F_1^T = Y,
\end{equation} 
where $F = F_1 \otimes F_1$, $\mathbf{x} = \vec(X)$, and $\mathbf{y} = \vec(Y)$. 
We furthermore assume that the regularization operator $R$ is defined by 
\begin{equation}
    R \mathbf{x} 
        = \begin{bmatrix} I \otimes R_1 \\ R_1 \otimes I \end{bmatrix} \vec(X) 
        = \begin{bmatrix} \vec(R_1X) \\ \vec(XR_1^T) \end{bmatrix},
\end{equation}
which corresponds to anisotropic regularization. 
Using some basic properties of the Kronecker product and the element-wise Hadamard product $\odot$, it can be shown that 
\begin{align} 
	F^T A F \mathbf{x} 
		& = \vec \left( F_1^T \left[ \tilde{A} \odot F_1 X F_1^T \right] F_1  \right), \label{eq:GD_eq1} \\ 
	R^T B R \mathbf{x} 
		& = \vec \left( \left[ \tilde{B}_1 \odot X R_1^T \right] R_1 \right) 
			+ \vec \left( R_1^T \left[ \tilde{B}_2 \odot R_1 X \right] \right), \label{eq:GD_eq2} \\ 
	\mathbf{b} 
		& = \vec \left( F_1^T \left[ \tilde{A} \odot X \right] F_1 \right), \label{eq:GD_eq3}
\end{align}
where $\tilde{A}$, $\tilde{B}_1$, and $\tilde{B}_2$ are such that $\vec( \tilde{A} ) = \boldsymbol{\alpha}$, $\vec(\tilde{B}_1) = \boldsymbol{\beta}^{1}$, and $\vec(\tilde{B}_2) = \boldsymbol{\beta}^{2}$, with $\boldsymbol{\beta} = [\boldsymbol{\beta}^{1},\boldsymbol{\beta}^{2}]$.
Combining \cref{eq:GD_eq1,eq:GD_eq2,eq:GD_eq3} yields 
\begin{equation}\label{eq:G_effic}
\resizebox{.9\textwidth}{!}{$\displaystyle 
	G \mathbf{x} = 
		\vec \left( F_1^T \left[ \tilde{A} \odot F_1 X F_1^T \right] F_1  \right) + 
		 \vec \left( \left[ \tilde{B}_1 \odot X R_1^T \right] R_1 \right) + 
		 \vec \left( R_1^T \left[ \tilde{B}_2 \odot R_1 X \right] \right) 
$}%
\end{equation} 
and therefore 
\begin{equation}\label{eq:grad_effic}
\begin{aligned} 
	\nabla J(\mathbf{x}) = 2 \Bigg[ 
	 	& \vec \left( F_1^T \left[ \tilde{A} \odot F_1 X F_1^T \right] F_1  \right) 
		+ \vec \left( \left[ \tilde{B}_1 \odot X R_1^T \right] R_1 \right) \\ 
		& + \vec \left( R_1^T \left[ \tilde{B}_2 \odot R_1 X \right] \right) 
		- \vec \left( F_1^T \left[ \tilde{A} \odot X \right] F_1 \right)
		\Bigg]. 
\end{aligned}
\end{equation} 
Observe that all of the matrices in \cref{eq:G_effic,eq:grad_effic} are significantly smaller than $F$ and $R$.

\subsection{Uncertainty quantification} 
\label{sub:UQ}

The proposed BCD algorithm has the advantage of allowing for uncertainty quantification in the reconstructed solution $\mathbf{x}$. 
For fixed $\boldsymbol{\alpha}$ and $\boldsymbol{\beta}$, Bayes' theorem and the conjugacy relationship \cref{eq:prod_likelihood_prior} yield 
\begin{equation}
    p(\mathbf{x}|\mathbf{y}) 
        \propto p(\mathbf{y}|\mathbf{x}) p(\mathbf{x}) 
        \propto \mathcal{N}(\mathbf{x}|\boldsymbol{\mu},C),  
\end{equation} 
where the mean $\boldsymbol{\mu}$ and the covariance matrix $C$ are again given by \cref{eq:mean_covariance2}. 
We can then sample from the normal distribution $\mathcal{N}(\boldsymbol{\mu},C)$ to obtain, for instance, credible intervals for every component of the solution $\mathbf{x}$. 
At the same time, we stress that this only allows for uncertainty quantification in $\mathbf{x}$ for given hyper-parameters $\boldsymbol{\alpha}$ and $\boldsymbol{\beta}$. 
The above approach does not include uncertainty in $\boldsymbol{\alpha}$ and $\boldsymbol{\beta}$ when these are treated as random variables themselves. 
This might be achieved by employing a computational more expensive sampling approach \cite{bardsley2012mcmc}, which we will investigate in future work.

\subsection{Relationship to current methodology}
\label{sec:othermethods}

We now address the connection between the proposed BCD algorithm and some existing methods.

\subsubsection{Iterative alternating sequential algorithm}

There are both notable similarities and key distinctions between the proposed BCD algorithm and the iterative alternating sequential (IAS) algorithm, developed in \cite{calvetti2007introduction,calvetti2007gaussian} and further investigated in \cite{calvetti2015hierarchical,calvetti2019hierachical}. 
Both algorithms estimate the unknown $\mathbf{x}$ and other involved parameters by alternatingly updating them. 
However, in contrast to the BCD method, the IAS algorithm assumes that the noise covariance matrix $A$ is known, which then allows the restriction to white Gaussian noise $\boldsymbol{\nu} \sim \mathcal{N}(\mathbf{0},I)$; see \cite[Section 2]{calvetti2019hierachical}. 
Moreover, the IAS algorithm builds upon a conditionally Gaussian prior for which the elements of the diagonal covariance matrix are gamma-distributed, rather than the elements of the diagonal \emph{inverse} covariance matrix as done here, which does not result in a conjugate hyper-prior. 
This makes the update steps for $\mathbf{x}$ and the hyper-parameters of the prior more complicated. 
Finally, the IAS algorithm solves for the MAP estimate of the posterior, which does not provide uncertainty quantification in the reconstructed solution. 
By contrast, the proposed BCD method grants access to the solution posterior $p(\mathbf{x}|\mathbf{y})$ for fixed hyper-parameters.

\subsubsection{Iteratively reweighted least squares} 

The update steps \cref{eq:update_x,eq:update_alpha,eq:update_beta} resulting from \cref{algo:BCD_mean} can be interpreted as an iteratively reweighted least squares (IRLS) algorithm \cite{daubechies2010iteratively}. 
The idea behind IRLS algorithms is to recover, for instance, a sparse solution by penalizing the components of $\mathbf{x}$ by weighting them individually and iteratively updating these weights. 
Indeed, the update steps \cref{eq:update_x,eq:update_alpha,eq:update_beta} resemble reweighted Tikhonov-regularization strategies. 
In this regard, the BCD method provides a solid Bayesian interpretation for commonly used reweighting choices and might be used to tailor these weights to specific statistical assumptions on the underlying problem.

\subsubsection{ARD/SBL optimization via iteratively re-weighted $\ell^1$-minimization} 

The first SBL algorithms used the same $\mathbf{x}$-update as in \cref{algo:BCD_mean}, but updated the noise and prior parameters $\alpha$, $\boldsymbol{\beta}$ using the evidence approach (expectation maximization) or the fixed-point approach, \cite{tipping2001sparse,mackay1992bayesian}. 
Although these methods can yield sparse solutions, they have no convergence guarantees and become prohibitively slow for large problems. 
Subsequently, in \cite{wipf2007new} it was demonstrated that the (type-II) evidence approach can be interpreted as a (type-I) MAP approach with a special non-factorable prior. 
With this insight in hand, a more efficient algorithm was then proposed to update $\mathbf{\beta}$ based on re-weighted $\ell^1$-minimization, which provably converges to a local maximum of the evidence $p( \mathbf{y} | \boldsymbol{\alpha}, \boldsymbol{\beta} )$ (see \cref{eq:evidence_marg}) with respect to $\boldsymbol{\beta}$. 
For the `almost' general regularization operators considered here, we cannot use the algorithm proposed in \cite{wipf2007new} since the evidence becomes improper if $\kernel(R) \neq \{\mathbf{0}\}$ (see \cref{app:improper_evidences}). 
By contrast, the $\alpha$- and $\boldsymbol{\beta}$-updates in \cref{algo:BCD_mean} are decoupled and based on respectively maximizing the full conditional posteriors \cref{eq:post2,eq:post3} (if we solve for the mode of the posterior $p(\mathbf{x}, \boldsymbol{\alpha}, \boldsymbol{\beta} | \mathbf{y})$) or computing the mean of the full conditional posteriors \cref{eq:post2,eq:post3} (if we solve for the mean of the posterior $p(\mathbf{x}, \boldsymbol{\alpha}, \boldsymbol{\beta} | \mathbf{y})$). 
We were able to derive explicit and efficient formulas for these based on the conditionally conjugate relationships between the likelihood, prior, and hyper-priors.   

\section{Numerical results} 
\label{sec:numerics} 

The MATLAB code used to generate the numerical tests presented here is open access and can be found at GitHub.\footnote{See \url{https://github.com/jglaubitz/generalizedSBL}}

\subsection{Computational complexity} 

We start with addressing the computational complexity of the proposed BCD algorithm (\cref{algo:BCD_mean}) for Bayesian inference. 
Assume that \cref{algo:BCD_mean} stops after $L$ iterations, either because the algorithm has converged or reached the maximum number of iterations. 
In every iteration, the algorithm performs the $\mathbf{x}$-update \cref{eq:update_x}, the $\mathbf{\alpha}$-update \cref{eq:update_alpha}, and the $\boldsymbol{\beta}$-update \cref{eq:update_beta}. 
Denoting their computational complexity by $\mathcal{O}(h_x)$, $\mathcal{O}(h_{\alpha})$, and $\mathcal{O}(h_{\beta})$, respectively, the total computational complexity of the BCD method is $\mathcal{O}(L ( h_x + h_{\alpha} + h_{\beta} ) )$.

\subsubsection*{The $\mathbf{x}$-update}

If $\mathbf{x} \in \R^n$ represents a one-dimensional signal and the $\mathbf{x}$-update \cref{eq:update_x} is solved using the PCD method, then the computational complexity of this update is $\mathcal{O}(\tilde{n})$, where $\tilde{n}$ is the number of the non-zero elements of the coefficient matrix $G \in \R^{n \times n}$ on the left-hand side of \cref{eq:update_x}.\footnote{This assumes that the coefficient matrix itself is computed in $\mathcal{O}(\tilde{n})$.}
On the other hand, if $\mathbf{x} = \Vec(X) \in \R^{n^2}$ is the vectorized representation of an image $X \in \R^{n \times n}$ and the coefficient matrix $G \in \R^{n^2 \times n^2}$ is dense.  In this case we solve the $\mathbf{x}$-updated \cref{eq:update_x} using the efficient gradient descent approach described in \cref{sub:efficient_xupdate}. 
This method has a computational complexity of $\mathcal{O}(n^3)$ for a fixed number of iterations.\footnote{In our implementation we used five gradient descent steps for each $\mathbf{x}$-update.} 
We thus have $h_x = \max\{n^3, \tilde{n}\}$.

\subsubsection*{The $\boldsymbol{\alpha}$- and $\boldsymbol{\beta}$-updates}

If $\mathbf{x} \in \R^n$, $F \in \R^{m \times n}$, and $R \in \R^{k \times n}$, then $\boldsymbol{\alpha}$, $\boldsymbol{\beta}$ in \cref{eq:update_alpha,eq:update_beta} can be computed in $\mathcal{O}(nm)$ and $\mathcal{O}(nk)$, respectively. 
Assuming that $F$ and $R$ only contain $\tilde{n}_F$ and $\tilde{n}_R$ elements, then the computational complexity of the $\boldsymbol{\alpha}$- and $\boldsymbol{\beta}$-updates reduces to $\mathcal{O}(\tilde{n}_F)$ and $\mathcal{O}(\tilde{n}_R)$, respectively.
We thus have $h_{\alpha} = \max\{ nm, \tilde{n}_F \}$ and $h_{\beta} = \max\{ nk, \tilde{n}_R \}$.

\subsection{Denoising a sparse signal} 

Consider the sparse nodal values $\mathbf{x}$ of a signal $x:[0,1] \to \R$ at $n=20$ equidistant points. 
All of the values in $\mathbf{x}$ are zero except at four randomly selected locations, where the values were set to $1$. 
We are given noisy observations $\mathbf{y}$ which result from adding i.\,i.\,d.\ zero-mean normal noise with variance $\sigma^2 = 5 \cdot 10^{-2}$ to the exact values $\mathbf{x}$. 
The signal-to-noise ratio (SNR), defined as $E[\mathbf{x}^2]/\sigma^2$ with $E[\mathbf{x}^2] = (x_1^2+\dots+x_n^2)/n$, is $4$. 

\begin{figure}[tb]
	\centering
  	\begin{subfigure}[b]{0.45\textwidth}
		\includegraphics[width=\textwidth]{%
      		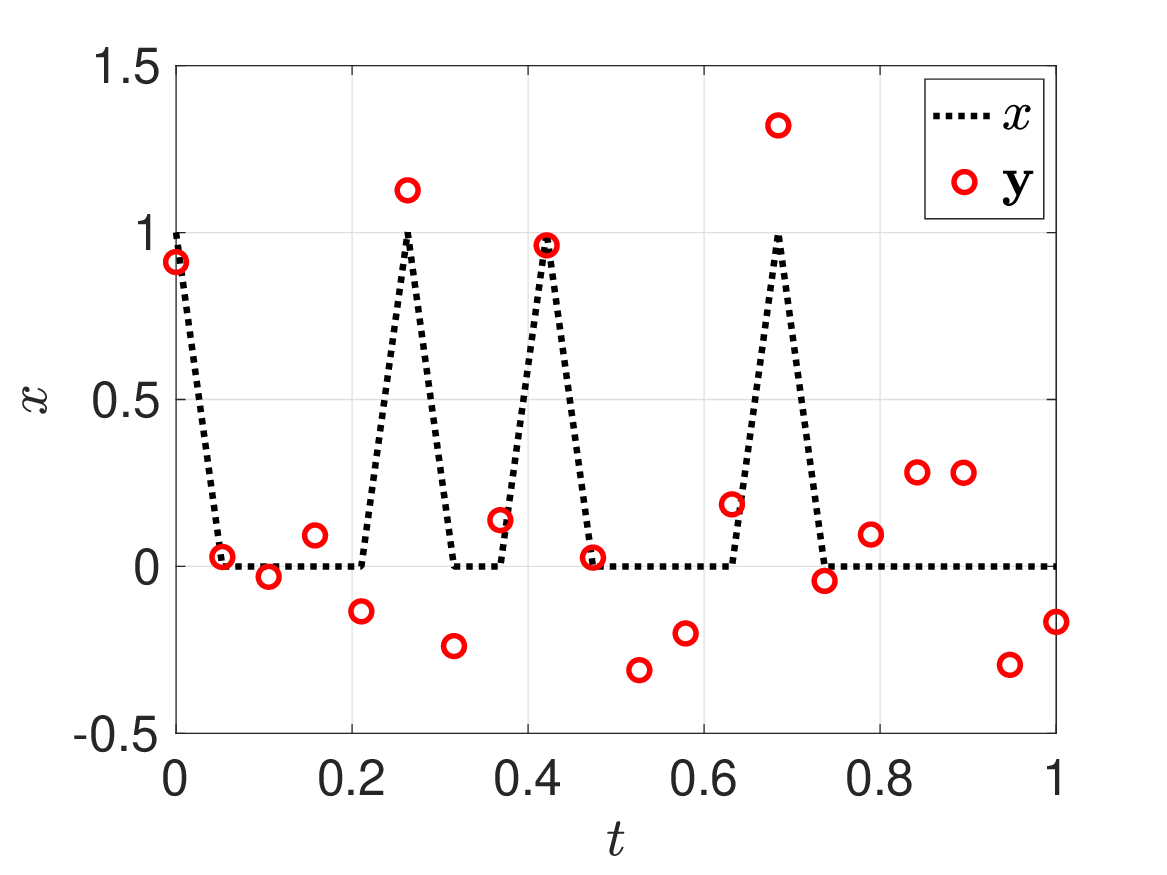} 
    	\caption{Signal $x$ and noisy observations $\mathbf{y}$}
    	\label{fig:denoising_1d_data}
  	\end{subfigure}%
  	~
  	\begin{subfigure}[b]{0.45\textwidth}
		\includegraphics[width=\textwidth]{%
      		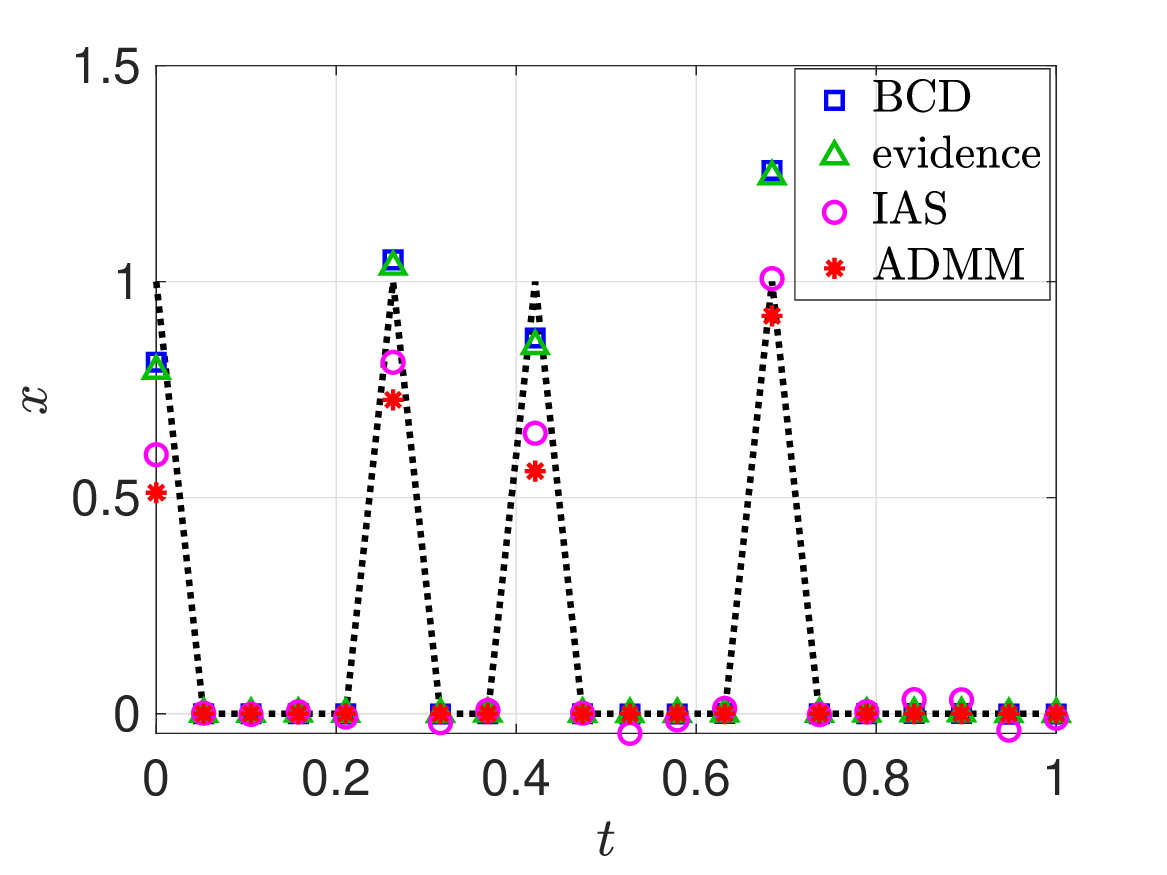} 
    		\caption{Reconstructions by different methods}
    		\label{fig:denoising_1d_reconstructions}
  	\end{subfigure}%
  	\caption{ 
  	The sparse signal $x$ and noisy observations $\mathbf{y}$ at $n=20$ equidistant points, and reconstructions by different methods 
  	}
  	\label{fig:denoising_1d}
\end{figure} 

\cref{fig:denoising_1d_data} illustrates the exact values of $x$ and the noisy observations $\mathbf{y}$. 
The corresponding data model and regularization operator are 
\begin{equation}
    \mathbf{y} = \mathbf{x} + \boldsymbol{\nu}, \quad R = I.
\end{equation} 
This simple test case allows us to compare the proposed BCD algorithm with some existing methods, some of which assume $\mathbf{x}$ itself to be sparse ($R = I$). 
\cref{fig:denoising_1d_reconstructions} provides a comparison of the BCD algorithm with (1) SBL using the evidence approach \cite{tipping2001sparse}, (2) the IAS method \cite{calvetti2007introduction,calvetti2007gaussian} solving for the MAP estimate of the posterior, and (3) the alternating direction method of multipliers (ADMM) \cite{boyd2011distributed} solving the deterministic $\ell^1$-regularized problem \cref{eq:l1_RIP}. 
The free parameters of the IAS algorithm were fine-tuned by hand and chosen as $\beta = 1.55$ and $\theta_j^* = 5 \cdot 10^{-2}$ for $j=1,\dots,n$; see \cite{calvetti2019hierachical} for more details on these parameters. 
The regularization parameter $\lambda$ in \cref{eq:l1_RIP} was also fine-tuned by hand and set to $\lambda = 2 \sigma^2 \|\mathbf{x}\|_0$. 
Finally, for the proposed BCD algorithm and the evidence approach, we assumed the noise variance $\sigma^2$ to be unknown, which therefore had to be estimated by the method as well. 
We can see in \cref{fig:denoising_1d_reconstructions} that for this example all of the SBL-based methods perform similarly. On the other hand, the ADMM yeilds a more regularized reconstruction, 
which might be explained by the uniform nature of the $\ell^1$-regularization term in \cref{eq:l1_RIP}. This is in contrast to the hierarchical Bayesian model which allows for spatially varying regularization. 
In this regard we note that there are weighted $\ell_1$-regularization methods \cite{candes2008enhancing,chartrand2008iteratively,adcock2019joint} that incorporate spatially varying regularization parameters. 
While such techniques can improve the resolution near the non-zero values in sparse signals, as well as near the internal edges in images, they are still point estimates and thus do not provide additional uncertainty information. 
Hence in the current investigation  we simply employ the standard ADMM with a fine-tuned non-varying regularization parameter as a reasonable comparison.

\subsection{Deconvolution of a piecewise constant signal}
\label{sub:deconvolution_1d}

\begin{figure}[tb]
	\centering
  	\begin{subfigure}[b]{0.45\textwidth}
		\includegraphics[width=\textwidth]{%
      		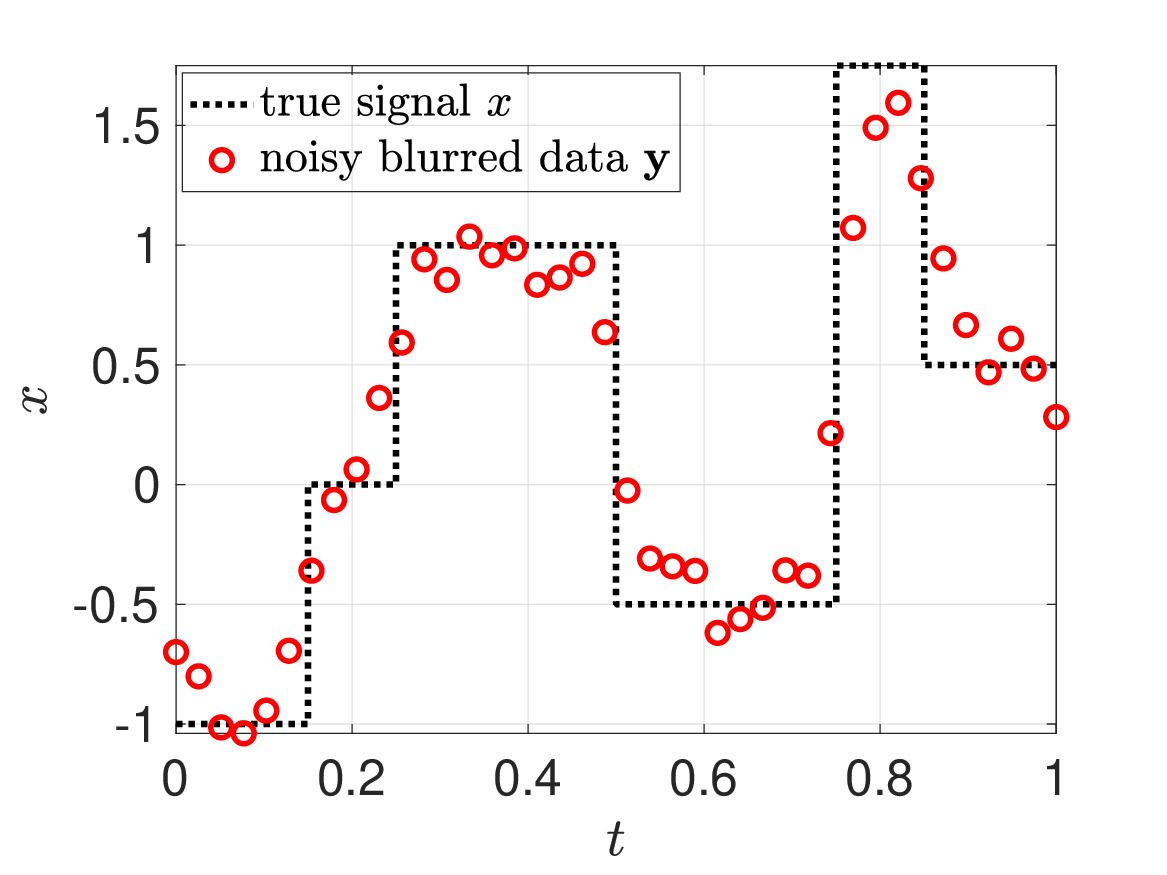} 
    	\caption{Signal $x$ and noisy blurred data $\mathbf{y}$}
    	\label{fig:deconvolution_1d_data_noise001}
  	\end{subfigure}%
  	~
  	\begin{subfigure}[b]{0.45\textwidth}
		\includegraphics[width=\textwidth]{%
      		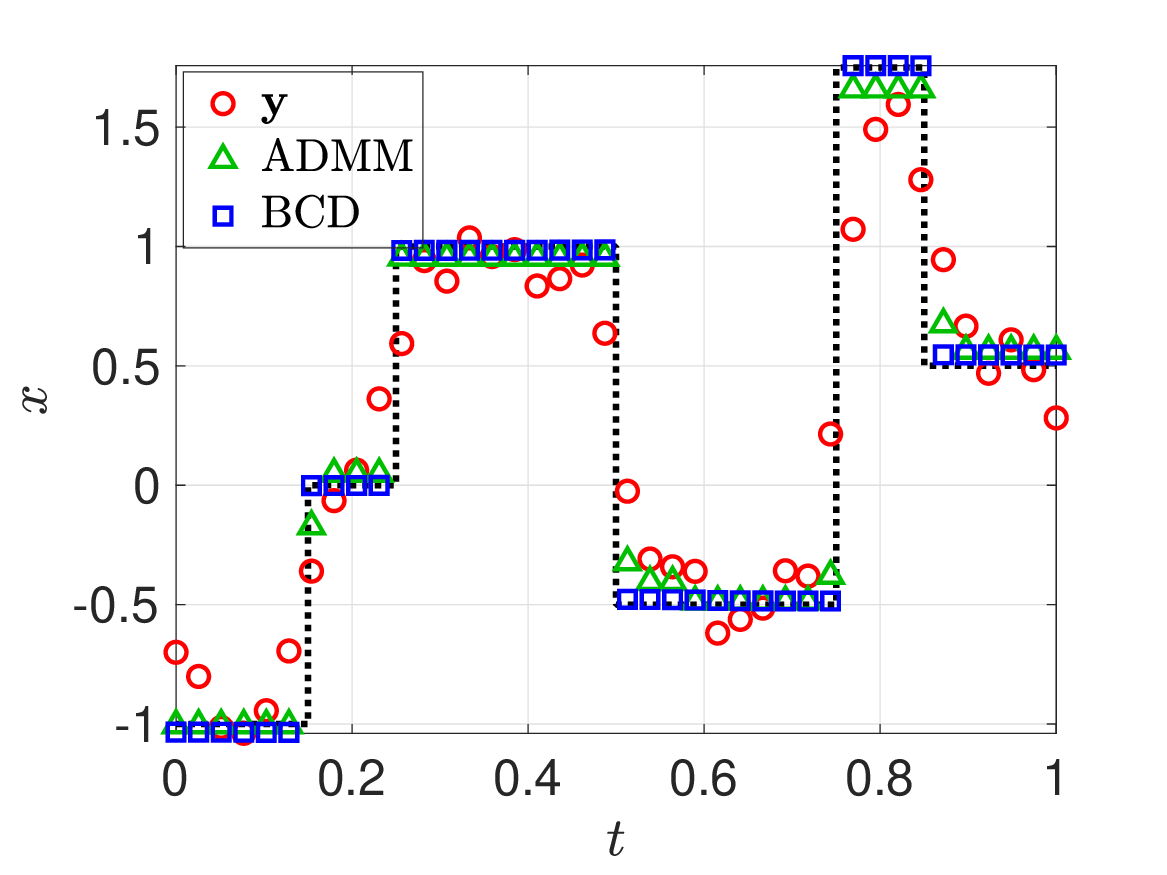} 
    		\caption{Different reconstructions}
    		\label{fig:deconvolution_1d_reconstructions_noise001}
  	\end{subfigure}%
	\\ 
	\begin{subfigure}[b]{0.45\textwidth}
		\includegraphics[width=\textwidth]{%
      		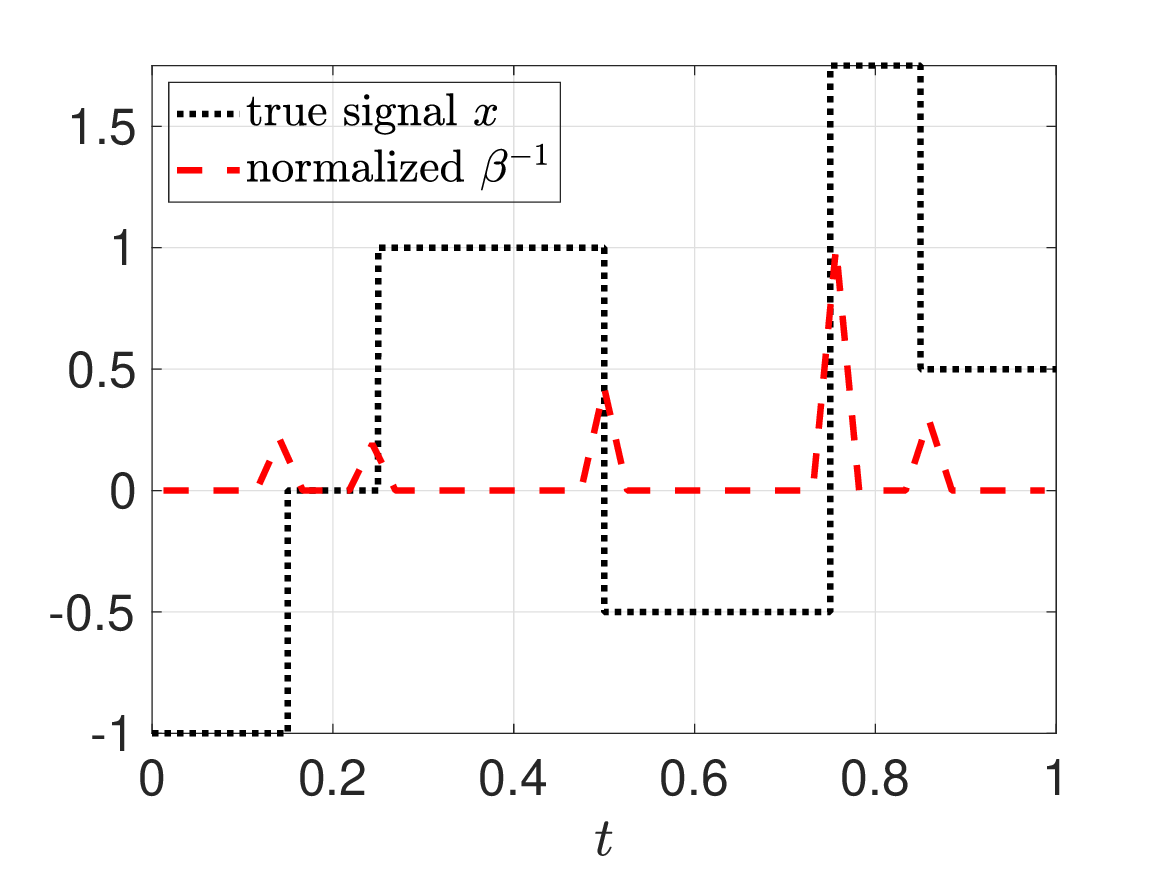} 
    	\caption{Normalized covariance parameter $\beta^{-1}$}
    	\label{fig:deconvolution_1d_beta_noise001}
  	\end{subfigure}%
  	~
  	\begin{subfigure}[b]{0.45\textwidth}
		\includegraphics[width=\textwidth]{%
      		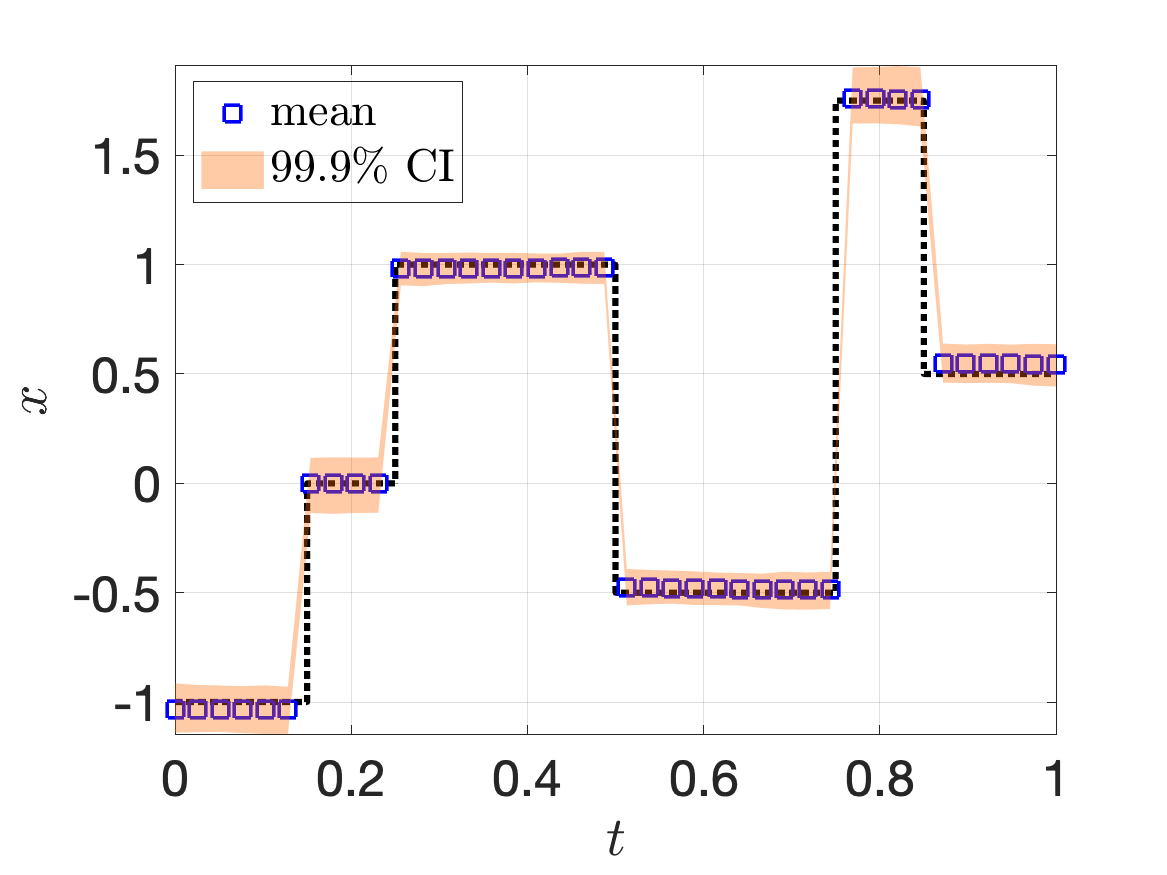} 
    		\caption{Mean and $99.9\%$ credible intervals}
    		\label{fig:deconvolution_1d_CI_noise001}
  	\end{subfigure}%
  	\caption{ 
  	Deconvolution of a piecewise constant signal $x$ from noisy blurred data $\mathbf{y}$ with i.\,i.\,d.\ zero-mean normal noise with variance $\sigma^2 = 10^{-2}$  
  	}
  	\label{fig:deconvolution_1d_noise001}
\end{figure}

We next consider deconvolution of the piecewise constant signal $x:[0,1] \to \R$ illustrated in \cref{fig:deconvolution_1d_noise001}. 
The corresponding data model and regularization operator are respectively given by 
\begin{equation}\label{eq:model_deconvolution_1d}
    \mathbf{y} = F \mathbf{x} + \boldsymbol{\nu}, \quad 
    R = 
    \begin{bmatrix}
        -1 & 1 & & \\ 
         & \ddots & \ddots & \\ 
         & & -1 & 1 
    \end{bmatrix} 
    \in \R^{(n-1) \times n},
\end{equation}
where $\boldsymbol{\nu} \sim \mathcal{N}(\mathbf{0},\sigma^2 I)$ with $\sigma^2 = 10^{-2}$ ($\mathrm{SNR} \approx 80$) and $F$ is obtained by applying the midpoint quadrature to the convolution equation 
\begin{equation} 
	y(s) = \int_0^1 k(s-s') x(s) \intd s'.
\end{equation}
We assume a Gaussian convolution kernel of the form 
\begin{equation} 
	k(s) = \frac{1}{2 \pi \gamma^2} \exp\left( - \frac{s^2}{2 \gamma^2} \right)
\end{equation} 
with blurring parameter $\gamma = 3 \cdot 10^{-2}$. 
The forward operator thus is 
\begin{equation}\label{eq:disc_convolution}
	[F]_{ij} = h k( h[i-j] ), \quad i,j=1,\dots,n,
\end{equation}
where $h=1/n$ is the distance between consecutive grid points. 
Note that $F$ has full rank but quickly becomes ill-conditioned. 

\cref{fig:deconvolution_1d_data_noise001} illustrates the true signal $x$ as well as the given noisy blurred data $\mathbf{y}$ at $n=40$ equidistant points. 
\cref{fig:deconvolution_1d_reconstructions_noise001} provides the reconstructions using the SBL-based BCD algorithm and the ADMM $\ell^1$-regularized inverse problem \cref{eq:l1_RIP}. 
The regularization parameter $\lambda$ in \cref{eq:l1_RIP} was again fine-tuned by hand and chosen as $\lambda = 2 \sigma^2 \|R\mathbf{x}\|_0$. 
We do not include any of the existing SBL algorithms considered before (the evidence approach and IAS algorithm) since they cannot be applied to the non-quadratic regularization operator $R$ in \cref{eq:model_deconvolution_1d} without modifying this operator first. 
\cref{fig:deconvolution_1d_beta_noise001} illustrates the normalized prior covariance parameters $\beta^{-1}$ which are estimated as part of the BCD algorithm. 
Observe that the values are significantly larger at the locations of the jump discontinuities. 
This allows the reconstruction to ``jump" and highlights the nonuniform character of regularization in the hierarchical Bayesian model suggested in \cref{sec:Bayesian}.  
Finally, \cref{fig:deconvolution_1d_CI_noise001} demonstrates the possibility to quantify uncertainty when using the BCD algorithm by providing the $99.9\%$ credible intervals of the solution posterior $p(\mathbf{x} | \mathbf{y})$ for the final estimates of $\alpha$ and $\boldsymbol{\beta}$. 
Note that these credible intervals, especially their width, indicate the amount of uncertainty in the reconstruction. 

\begin{figure}[tb]
	\centering
  	\begin{subfigure}[b]{0.45\textwidth}
		\includegraphics[width=\textwidth]{%
      		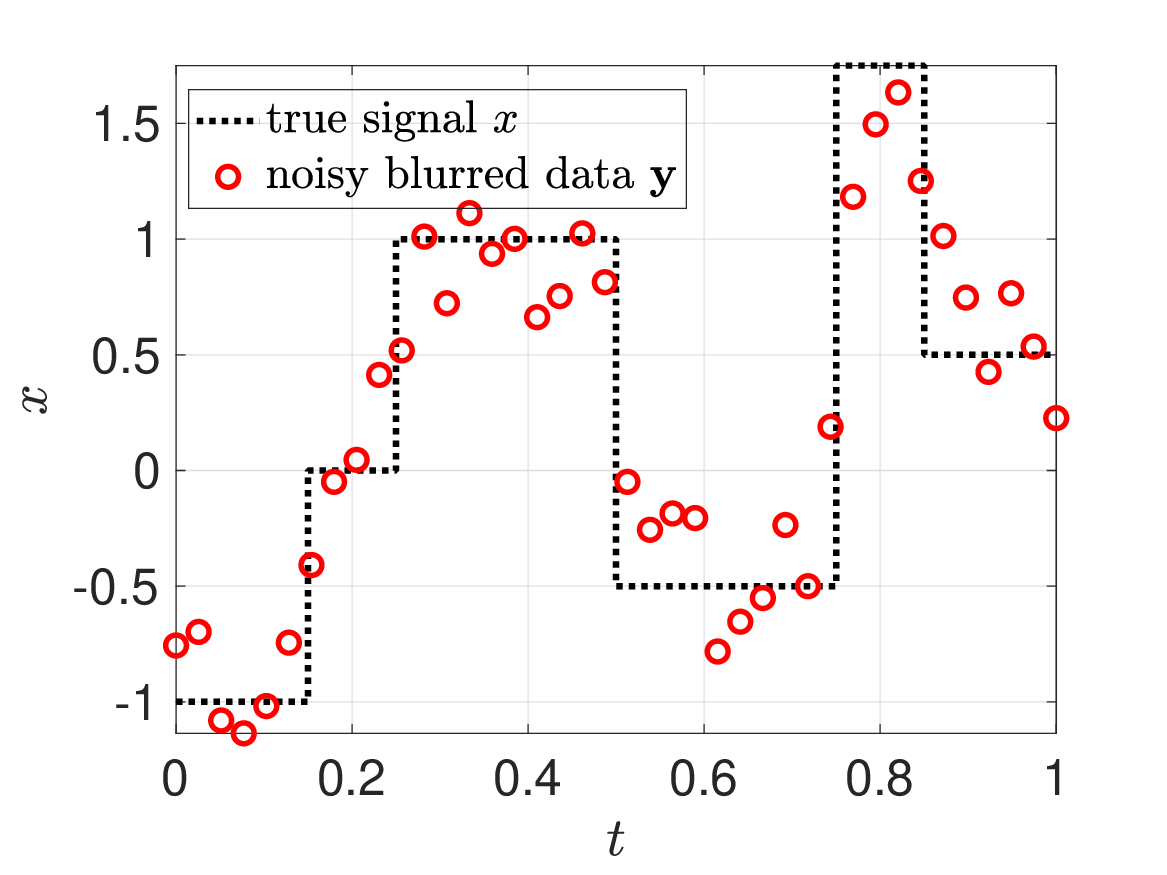} 
    	\caption{Signal $x$ and noisy blurred data $\mathbf{y}$}
    	\label{fig:deconvolution_1d_data_noise005}
  	\end{subfigure}%
  	~
  	\begin{subfigure}[b]{0.45\textwidth}
		\includegraphics[width=\textwidth]{%
      		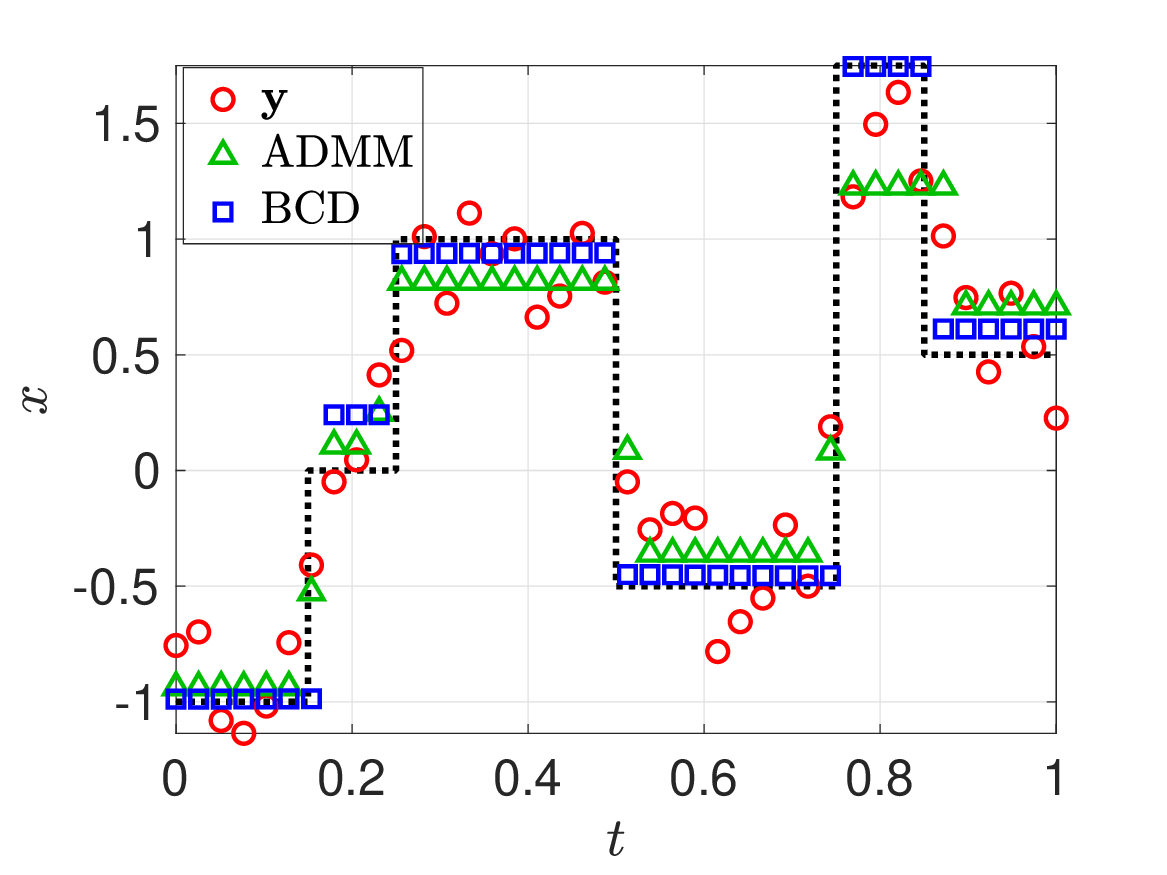} 
    		\caption{Different reconstructions}
    		\label{fig:deconvolution_1d_reconstructions_noise005}
  	\end{subfigure}%
	\\ 
	\begin{subfigure}[b]{0.45\textwidth}
		\includegraphics[width=\textwidth]{%
      		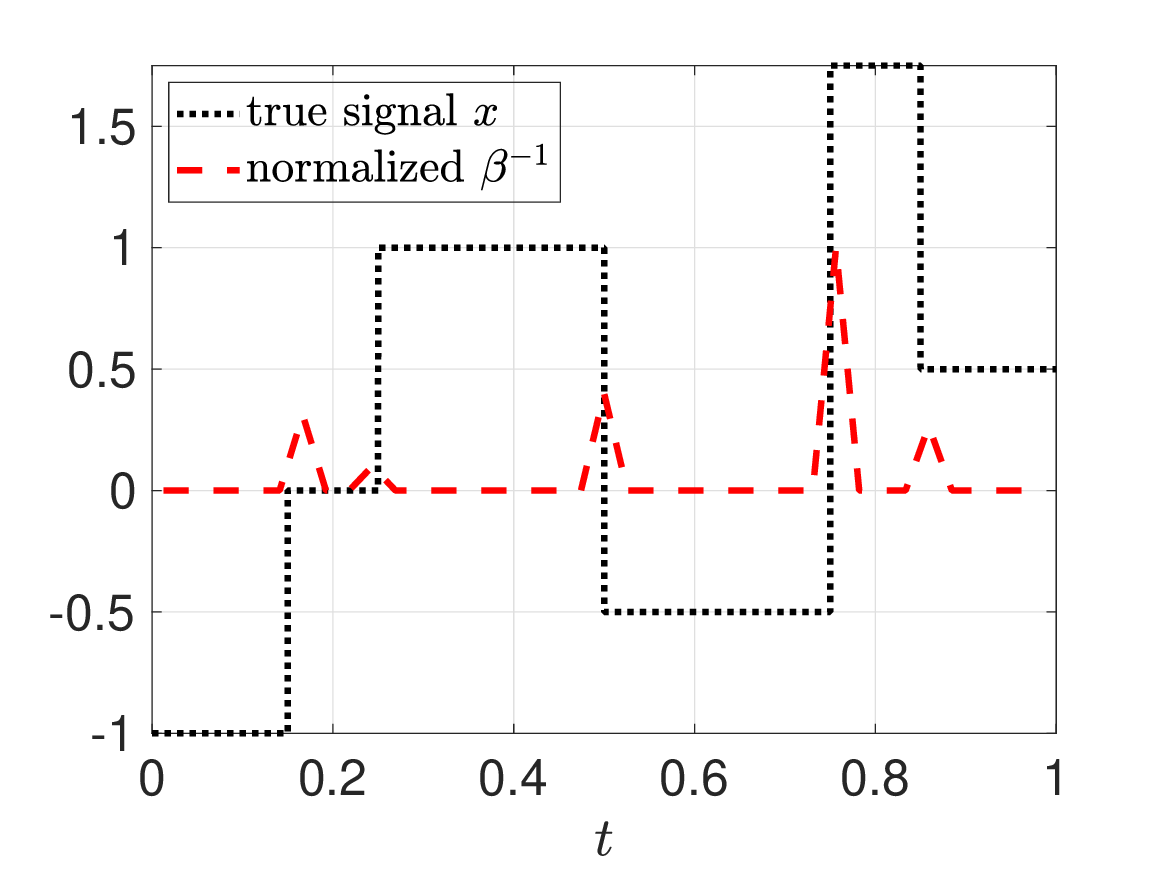} 
    	\caption{Normalized covariance parameter $\beta^{-1}$}
    	\label{fig:deconvolution_1d_beta_noise005}
  	\end{subfigure}%
  	~
  	\begin{subfigure}[b]{0.45\textwidth}
		\includegraphics[width=\textwidth]{%
      		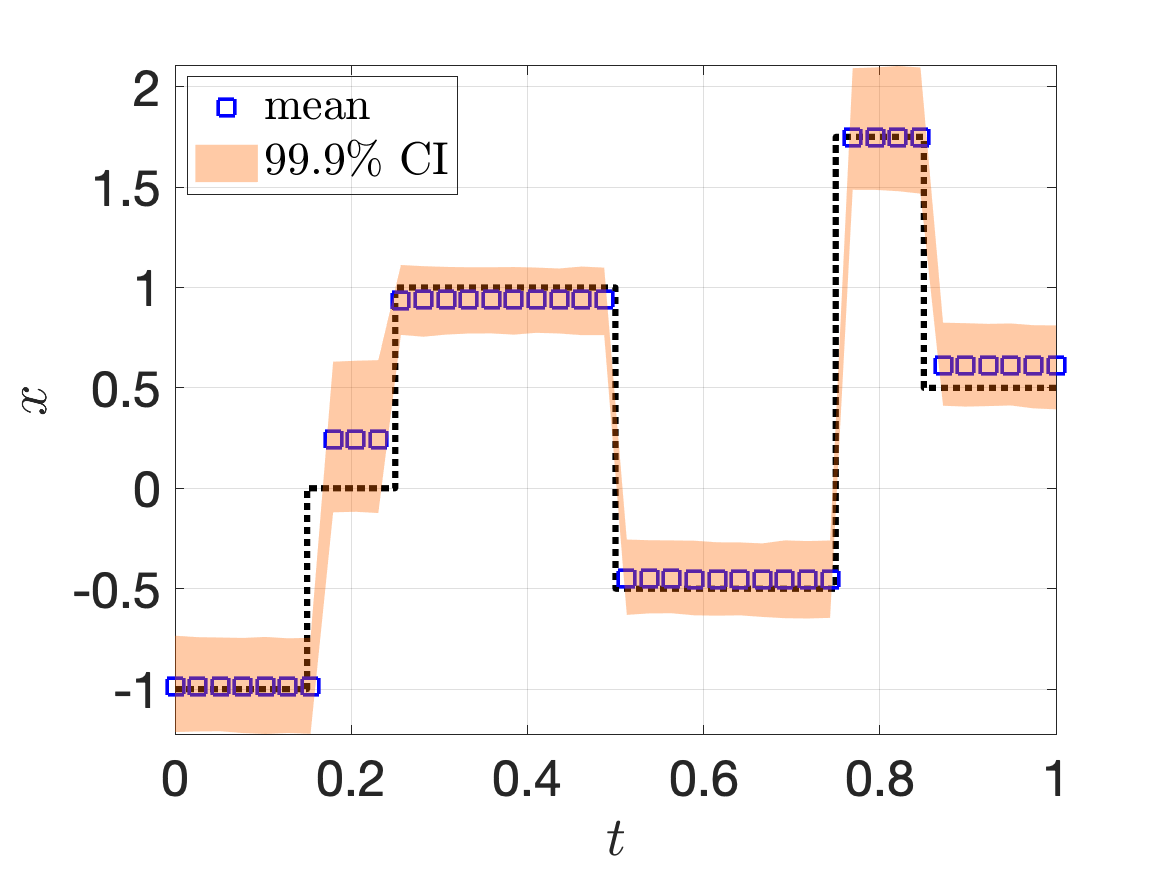} 
    		\caption{Mean and $99.9\%$ credible intervals}
    		\label{fig:deconvolution_1d_CI_noise005}
  	\end{subfigure}%
  	\caption{ 
  	Deconvolution of a piecewise constant signal $x$ from noisy blurred data $\mathbf{y}$ with i.\,i.\,d.\ zero-mean normal noise with variance $\sigma^2 = 5 \cdot 10^{-2}$ 
  	}
  	\label{fig:deconvolution_1d_noise005}
\end{figure} 

The results displayed in \cref{fig:deconvolution_1d_noise005} are for the same model with the noise variance increased by $500\%$, to $\sigma^2 = 5 \cdot 10^{-2}$ ($\mathrm{SNR} \approx 16$). 
The BCD algorithm now yields a less accurate reconstruction, especially between $t=0.15$ and $t=0.25$. 
This is also reflected in the corresponding normalized prior covariance parameters $\beta^{-1}$, which can be found \cref{fig:deconvolution_1d_beta_noise005}. 
Observe that the second peak around $t=0.25$ is underestimated and therefore causes the block associated with the region $[0.15,0.25]$ to be drawn towards the subsequent block associated with the region $[0.25,0.5]$. 
The increased uncertainty of the reconstruction is indicated by the $99.9\%$ credible intervals in \cref{fig:deconvolution_1d_CI_noise005}. 
In particular, we note the increased width of the credible interval in the region $[0.15,0.25]$.

\subsection{Combining different regularization operators} 

\begin{figure}[tb]
	\centering
  	\begin{subfigure}[b]{0.45\textwidth}
		\includegraphics[width=\textwidth]{%
      		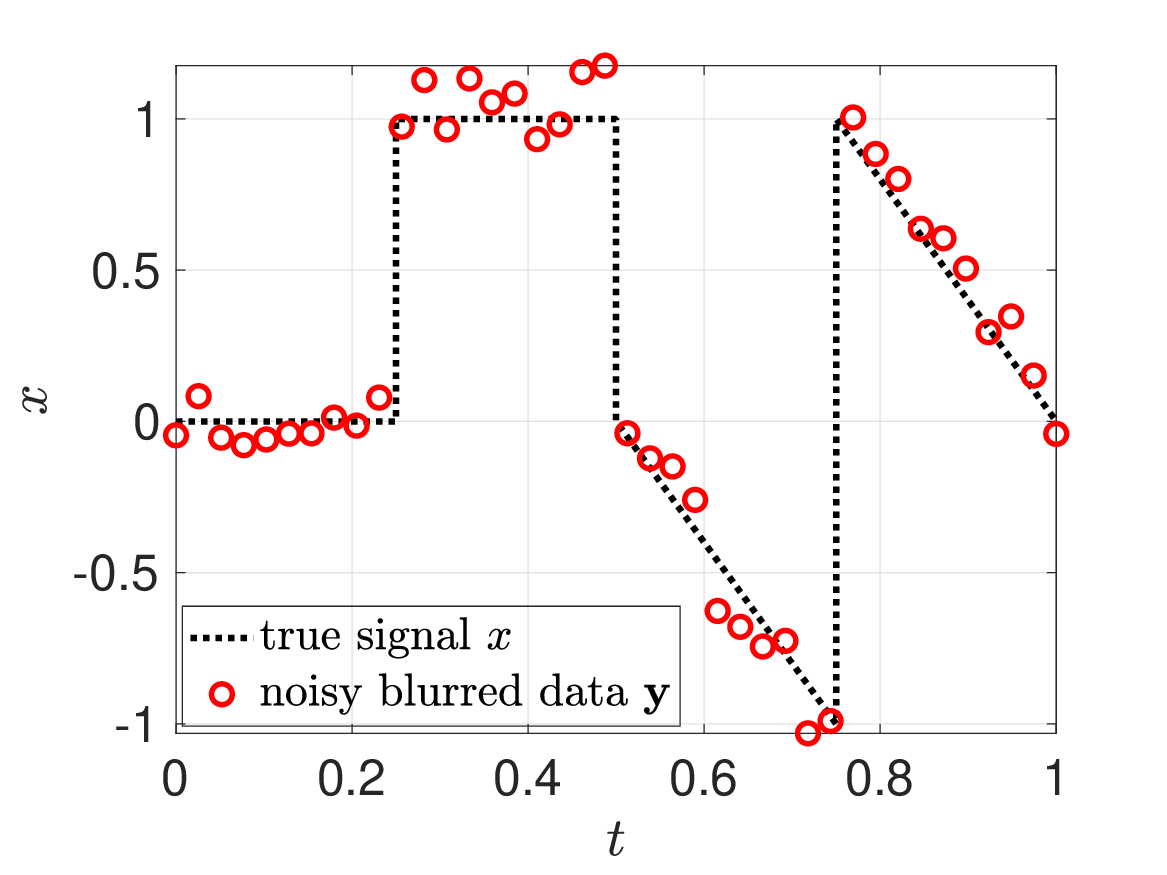} 
    	\caption{Signal $x$ and noisy observations $\mathbf{y}$}
    	\label{fig:deconvolution_combi_data}
  	\end{subfigure}%
  	~
  	\begin{subfigure}[b]{0.45\textwidth}
		\includegraphics[width=\textwidth]{%
      		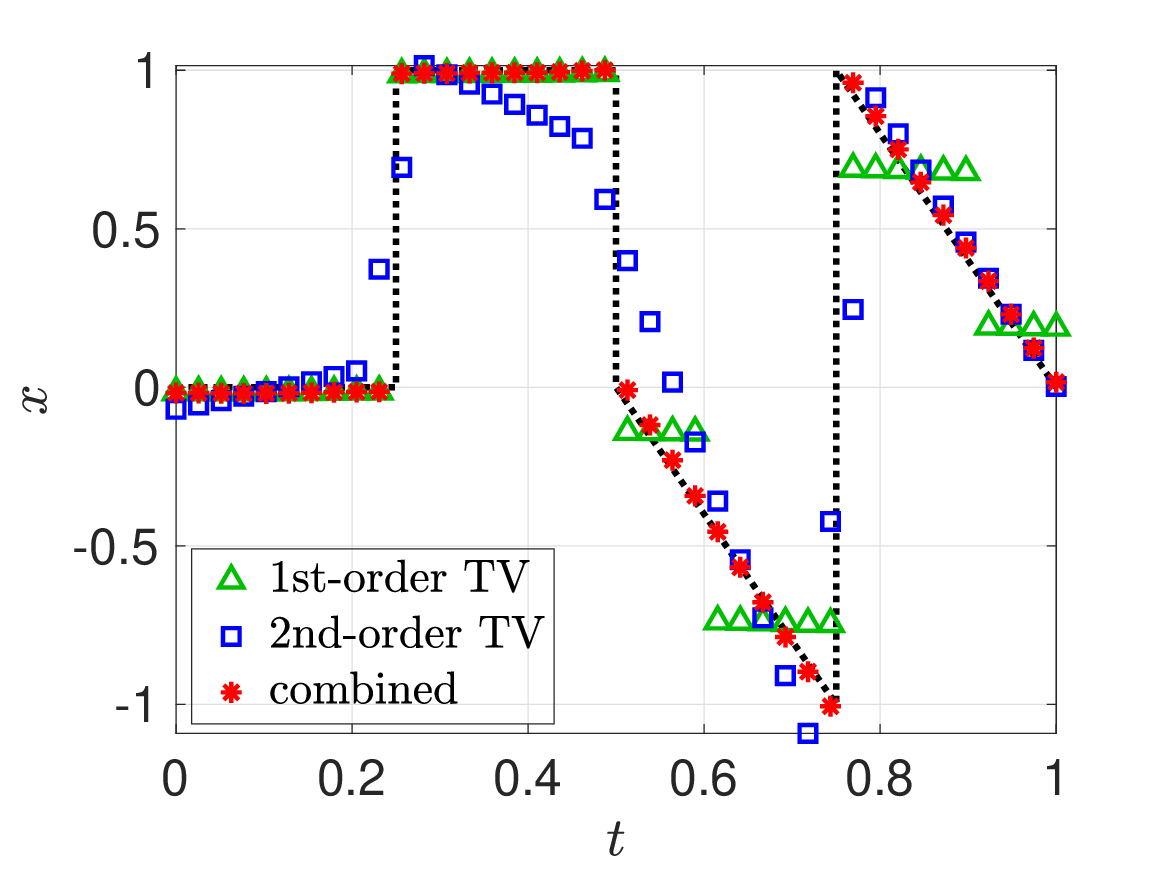} 
    		\caption{Reconstructions using different regualrizations}
    		\label{fig:deconvolution_combi_reconstructions}
  	\end{subfigure}%
  	\caption{ 
  	Signal $x$ and noisy observations $\mathbf{y}$ at $n=20$ equidistant points, and reconstructions by different methods 
  	}
  	\label{fig:deconvolution_combi}
\end{figure}

We next demonstrate that generalized SBL allows us to consider combinations of different regularization operators. 
Consider the signal $x:[0,1] \to \R$ illustrated in \cref{fig:deconvolution_combi_data}, which is piecewise constant on $[0,0.5]$ and piecewise linear on $[0.5,1]$. 
The corresponding data model is the same as before with convolution parameter $\gamma = 10^{-2}$ and i.\,i.\,d.\ zero-mean normal noise with variance $\sigma^2 = 10^{-2}$ ($\mathrm{SNR} \approx 40$). 
\cref{fig:deconvolution_combi_reconstructions} illustrates the reconstructions obtained by the BCD algorithm using a first- and second-order TV-regualrization operator, 
\begin{equation}\label{eq:reg_ops}
    R_1 = 
    \begin{bmatrix}
        -1 & 1 & & \\ 
         & \ddots & \ddots & \\ 
         & & -1 & 1 
    \end{bmatrix}, \quad 
    R_2 = 
    \begin{bmatrix}
        -1 & 2 & -1 & & \\ 
         & \ddots & \ddots & \ddots & \\ 
         & & -1 & 2 & -1 
    \end{bmatrix}, 
\end{equation}
which promote piecewise constant and piecewise linear solutions, respectively. 
Observe that neither $R_1$ nor $R_2$ is even square, meaning that both would have to be modified by introducing additional rows to apply a standard SBL approach, which can become increasingly complicated for higher orders and multiple dimensions. 

It is evident from \cref{fig:deconvolution_combi_reconstructions} that using first-order TV-regularization yields a less accurate reconstruction in $[0.5,1]$, where the signal is piecewise linear,\footnote{This well-known artifact of first-order TV-regularization is often called the ``staircasing'' effect and motivates using higher order TV-regularization, \cite{archibald2016image,stefan2010improved}.} while using second-order TV-regularization yields a less accurate reconstruction in $[0,0.5]$, where the signal is piecewise constant. 
However, generalized SBL and the proposed BCD algorithm allows us to consider the combined regularization operator 
\begin{equation}\label{eq:reg_op2}
    R = 
    \begin{bmatrix}
        -1 & 1 & & & & & \\ 
         & \ddots & \ddots & & & & \\ 
         & & -1 & 1 & & & \\ 
         & & -1 & 2 & -1 & & \\ 
         & & & \ddots & \ddots & \ddots & \\ 
         & & & & -1 & 2 & -1 \\ 
    \end{bmatrix} 
    \in \R^{(n-3) \times n}.
\end{equation} 
Assuming $n = 2q$, the first $k-1$ rows correspond to first-order TV-regularization while the last $k-2$ rows correspond to second-order TV-regularization. 
The advantage of using this nonstandard regularization operator in the BCD algorithm is demonstrated by the red stars in \cref{fig:deconvolution_combi_reconstructions}.

\subsection{Image deconvolution} 

\begin{figure}[tb]
	\centering
  	\begin{subfigure}[b]{0.45\textwidth}
		\includegraphics[width=\textwidth]{%
      		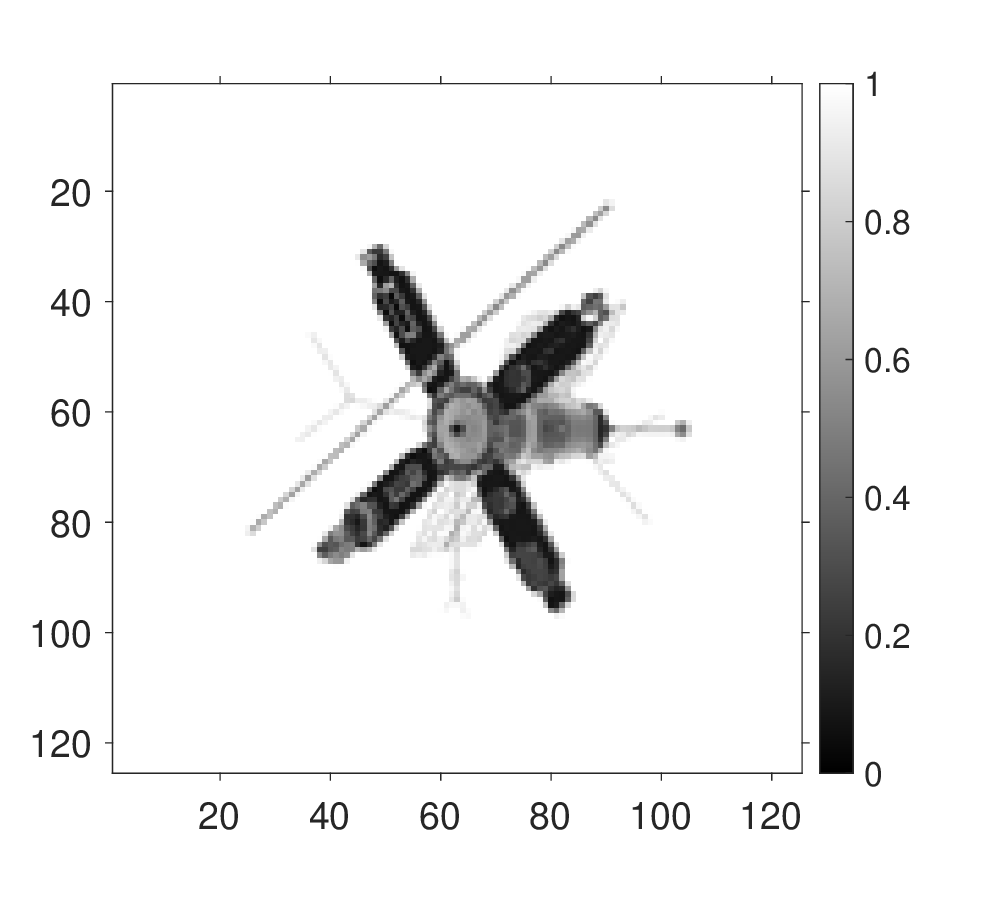} 
    	\caption{Reference image}
    	\label{fig:deconvolution_2d_ref}
  	\end{subfigure}%
  	~
  	\begin{subfigure}[b]{0.45\textwidth}
		\includegraphics[width=\textwidth]{%
      		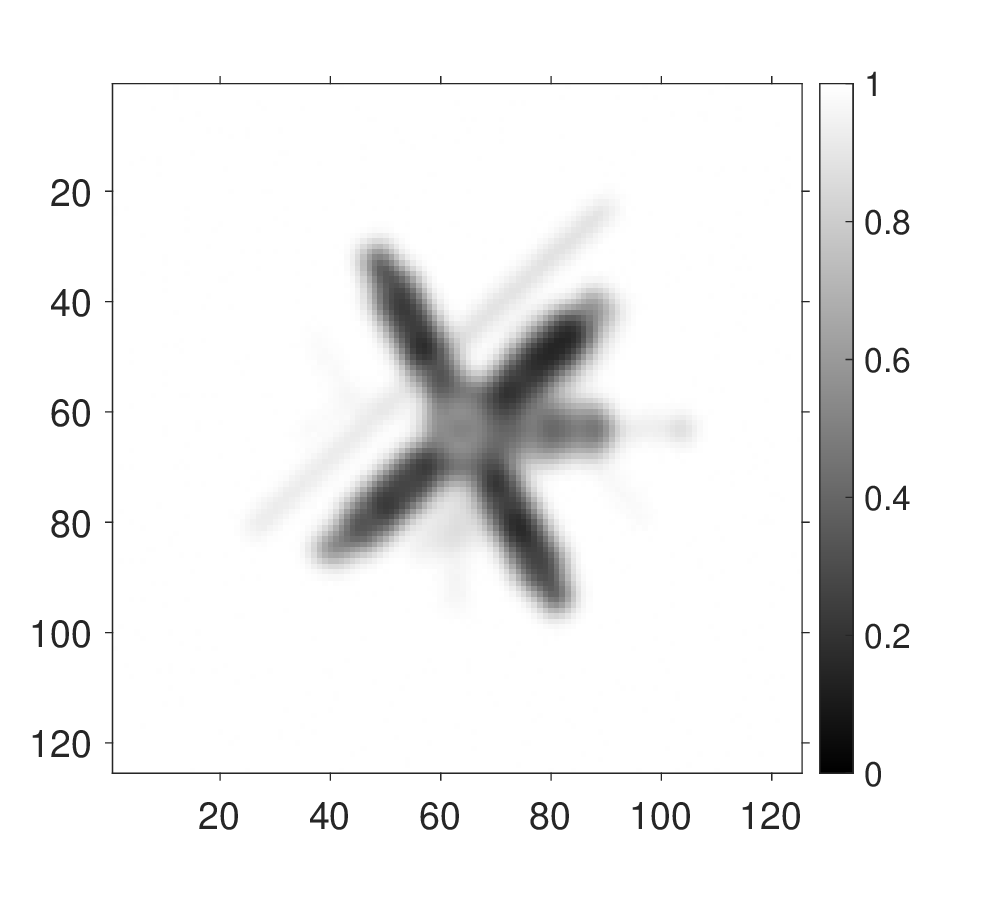} 
    		\caption{Noisy blurred image}
    		\label{fig:deconvolution_2d_blurred}
  	\end{subfigure}%
	\\ 
	\begin{subfigure}[b]{0.45\textwidth}
		\includegraphics[width=\textwidth]{%
      		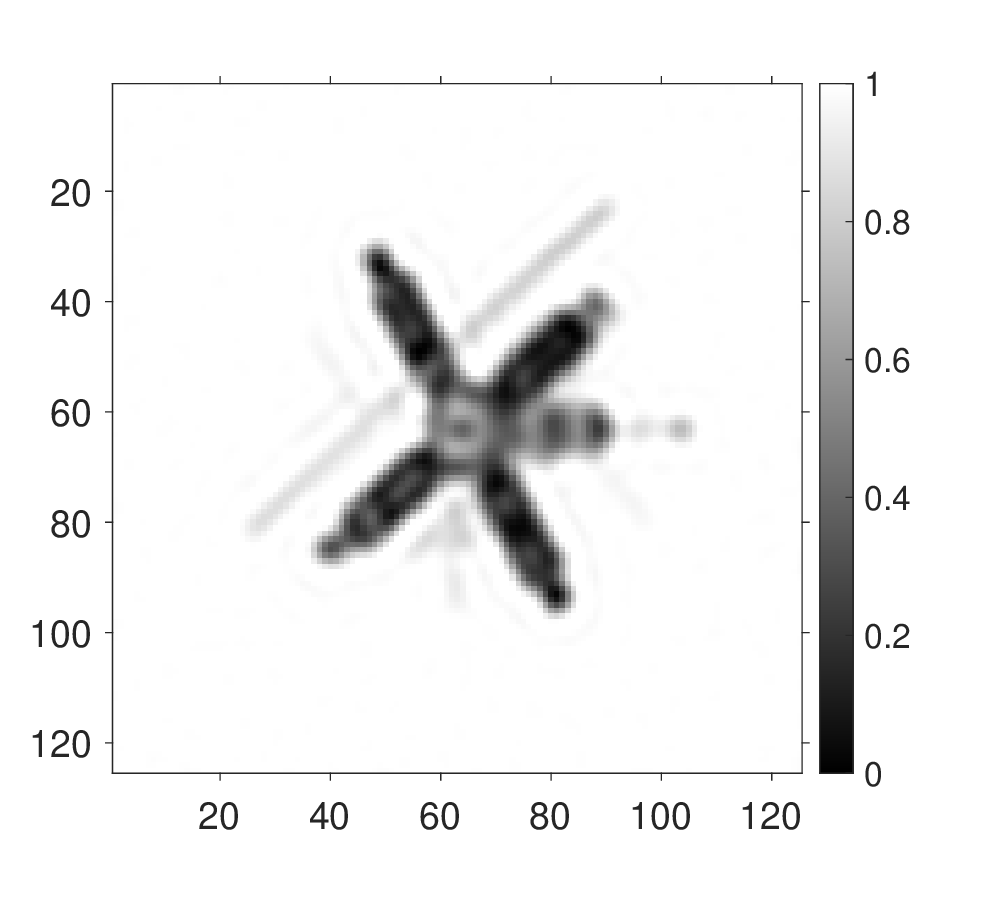} 
    	\caption{$\ell^1$-regularization by ADMM}
    	\label{fig:deconvolution_2d_ADMM}
  	\end{subfigure}%
  	~
  	\begin{subfigure}[b]{0.45\textwidth}
		\includegraphics[width=\textwidth]{%
      		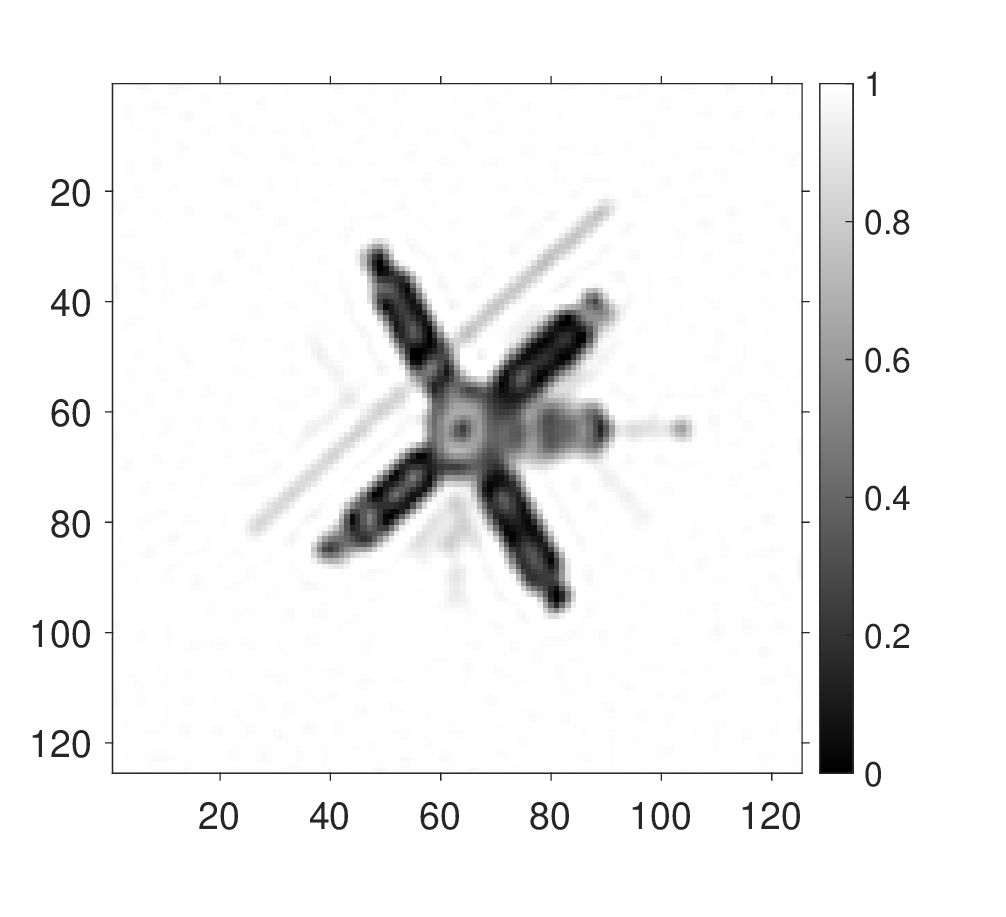} 
    		\caption{SBL by the BCD algorithm}
    		\label{fig:deconvolution_2d_BCD}
  	\end{subfigure}%
  	\caption{ 
  	The reference image, the corresponding noisy blurred image, and reconstructions using the ADMM and the BCD algorithm \cref{algo:BCD_mean} 
  	}
  	\label{fig:deconvolution_2d}
\end{figure} 

We next consider the reference image $X$ in \cref{fig:deconvolution_2d_ref} and its noisy blurred version $Y$ in \cref{fig:deconvolution_2d_blurred}. 
$Y$ results from $X$ by applying the discrete one-dimensional convolution operator \cref{eq:disc_convolution} in the two canonical coordinate directions and then adding i.\,i.\,d.\ zero-mean normal noise. 
The corresponding forward model is $Y = F X F^T + N$ or, equivalently, 
\begin{equation}
    \mathbf{y} = G \mathbf{x} + \boldsymbol{\nu}, 
\end{equation} 
after vectorization. 
Here, $\mathbf{z} = \Vec(Z)$ denotes the $mn\times1$ column vectors obtained by stacking the columns of the $m \times n$ matrix $Z$ on top of one another, and $G = F \otimes F$.  
Further, the blurring parameter and noise variance were chosen as $\gamma = 1.5 \cdot 10^{-2}$ and $\sigma^2 = 10^{-5}$ ($\mathrm{SNR} \approx 4 \cdot 10^3$) to make the test case comparable to the one in \cite[Section 4.2]{bardsley2012mcmc}.

\cref{fig:deconvolution_2d_ADMM,fig:deconvolution_2d_BCD} show the reconstructions obtained by the ADMM applied to \cref{eq:l1_RIP} and the SBL-based BCD algorithm with an anisotropic second-order TV operator 
\begin{equation}
    R = \begin{bmatrix} I \otimes D \\ D \otimes I \end{bmatrix} 
    \quad \text{with} \quad 
    D = 
    \begin{bmatrix}
        -1 & 2 & -1 & & \\ 
         & \ddots & \ddots & \ddots & \\ 
         & & -1 & 2 & -1
    \end{bmatrix} 
    \in \R^{(n-2) \times n}. 
\end{equation}
The regularization parameter $\lambda$ in \cref{eq:l1_RIP} was again fine-tuned by hand and set to $\lambda = 10^{-5}$. 
The BCD algorithm provides a sharper reconstruction (see \cref{fig:deconvolution_2d}) than the ADMM applied to the $\ell^1$-regularized inverse problem \cref{eq:l1_RIP}. 
Further parameter tuning might increase the accuracy of the reconstruction by the ADMM. 
By contrast, it is important to stress that the BCD algorithm requires no such exhaustive parameter tuning.

\subsection{Noisy and incomplete Fourier data} 
\label{sec:fourdata}

\begin{figure}[tb]
	\centering
  	\begin{subfigure}[b]{0.45\textwidth}
		\includegraphics[width=\textwidth]{%
      		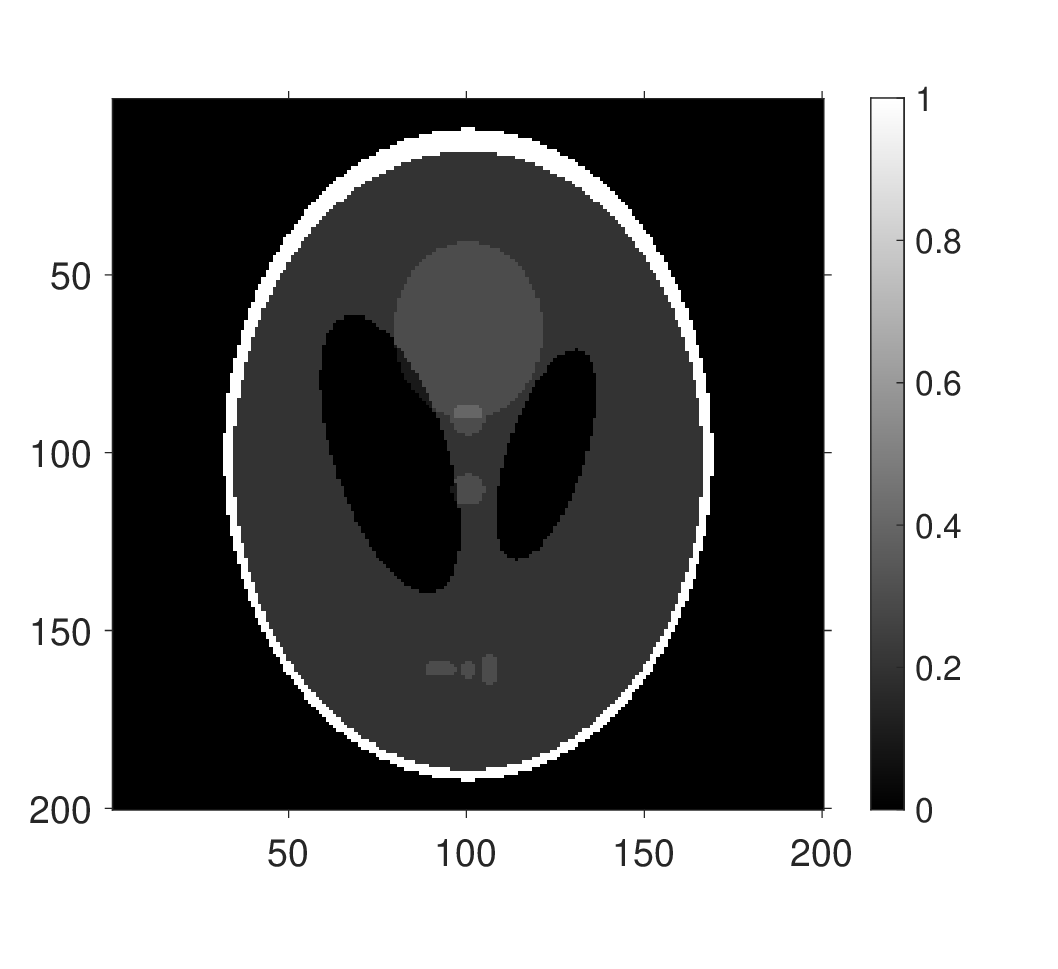} 
    	\caption{Reference image}
    	\label{fig:MRI_2d_ref}
  	\end{subfigure}%
  	~
  	\begin{subfigure}[b]{0.45\textwidth}
		\includegraphics[width=\textwidth]{%
      		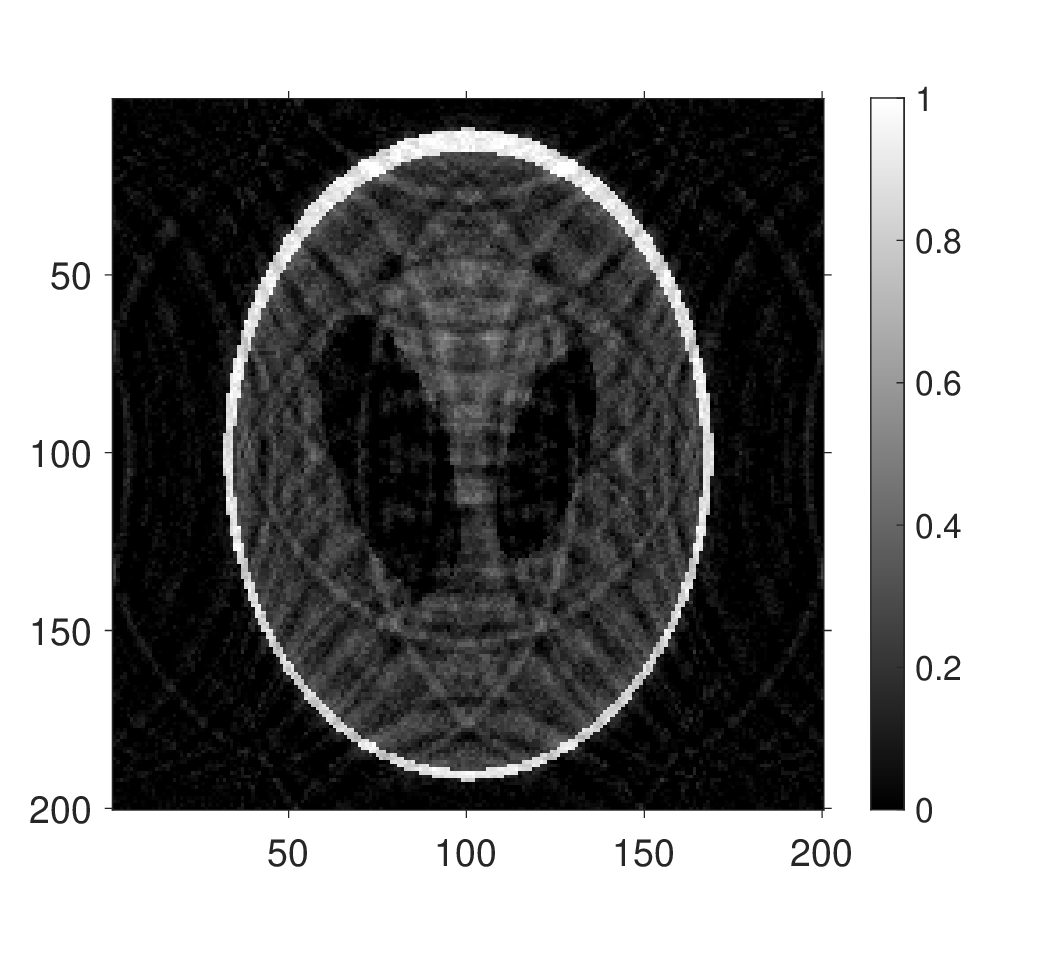} 
    		\caption{ML/LS estimate}
    		\label{fig:MRI_2d_ML}
  	\end{subfigure}%
	\\ 
	\begin{subfigure}[b]{0.45\textwidth}
		\includegraphics[width=\textwidth]{%
      		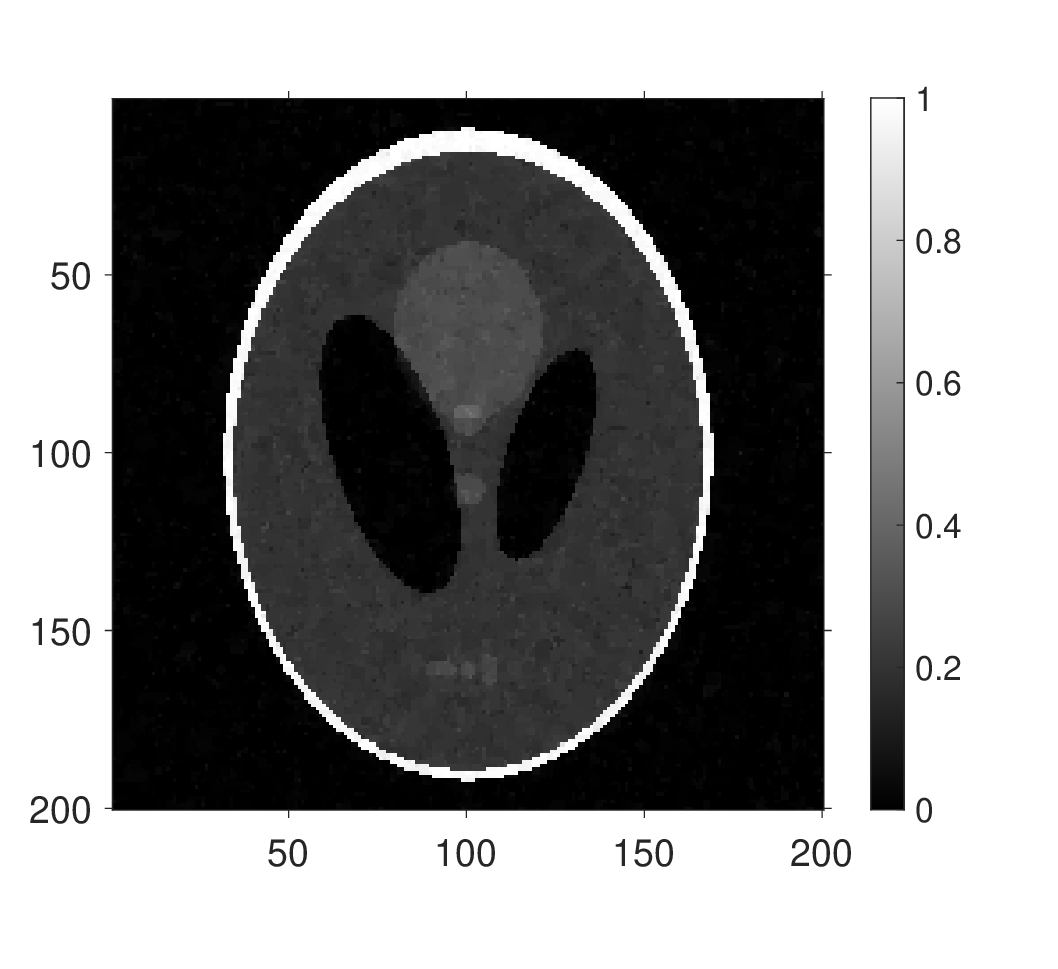} 
    	\caption{$\ell^1$-regularization by ADMM}
    	\label{fig:MRI_2d_ADMM}
  	\end{subfigure}%
  	~
  	\begin{subfigure}[b]{0.45\textwidth}
		\includegraphics[width=\textwidth]{%
      		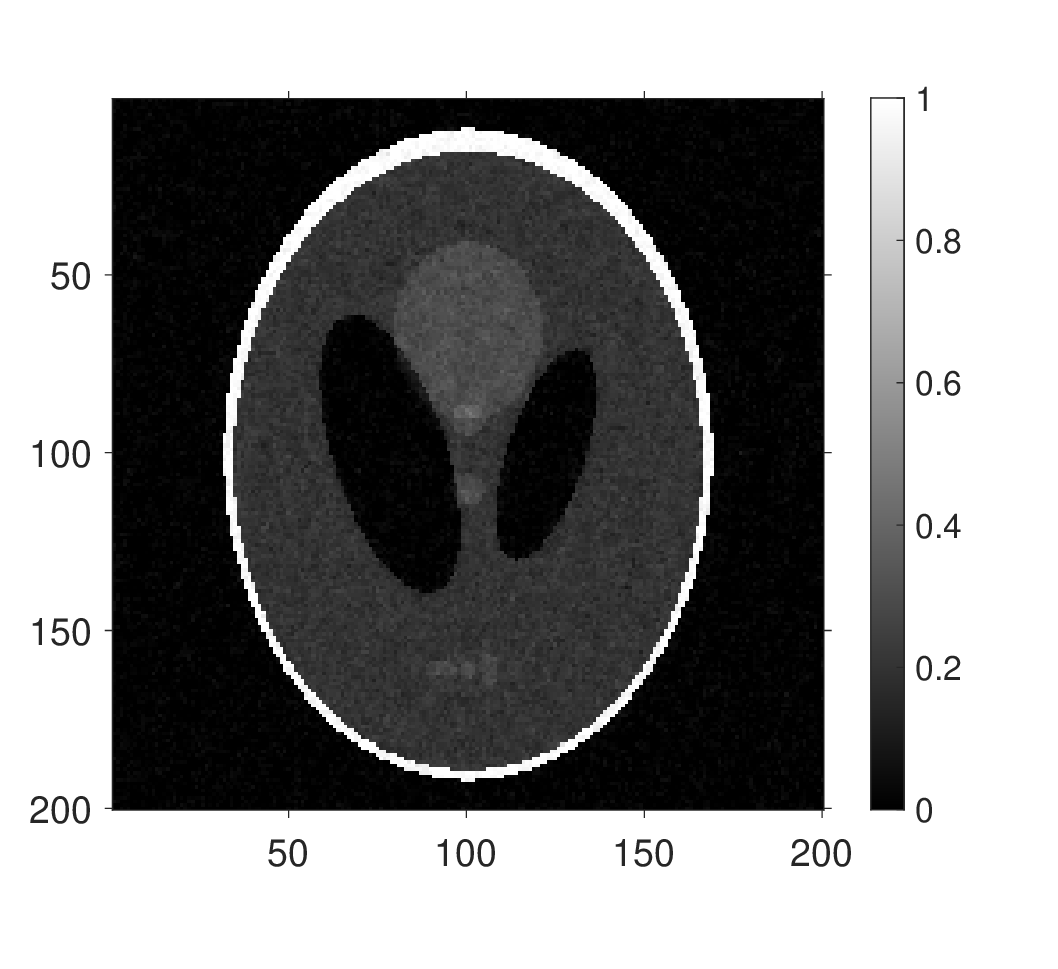} 
    		\caption{SBL by the BCD algorithm}
    		\label{fig:MRI_2d_BCD}
  	\end{subfigure}%
  	\caption{ 
  	(a) The Shepp--Logan phantom test image; (b) the ML/LS estimate, and reconstructions using (c) the ADMM applied to \cref{eq:l1_RIP}; and (d) the SBL-based BCD algorithm 
  	}
  	\label{fig:MRI_2d}
\end{figure}

We next address the reconstruction of images based on noisy and incomplete Fourier data, which is common in applications such as magnetic resonance imaging (MRI) and synthetic aperture radar (SAR). 
The popular prototype Shepp--Logan phantom test image is displayed in \cref{fig:MRI_2d_ref}.

The indirect data $\mathbf{y} = \Vec(Y)$ is given by applying the two-dimensional discrete Fourier transform to the reference image $X$, removing certain frequencies, and adding noise. 
Since in this investigation we are assuming $\mathbf{x} \in \mathbb{R}^n$, we consider the data model 
\begin{equation}\label{eq:data_model_Fourier}
    \begin{bmatrix} \mathrm{Re}(\mathbf{y}) \\ \mathrm{Im}(\mathbf{y}) \end{bmatrix} 
    = 
    \begin{bmatrix} \mathrm{Re}(G) \\ \mathrm{Im}(G) \end{bmatrix} 
    \mathbf{x} + \boldsymbol{\nu} 
\end{equation} 
with $\mathrm{Re}(\mathbf{y})$ and $\mathrm{Im}(\mathbf{y})$ respectively denoting the real and imaginary part of $\mathbf{y} \in \C^m$.\footnote{Our technique is not limited to real-valued solutions, and we will consider complex-valued solutions, such as those occurring in SAR, in future work.}  
Further, $\boldsymbol{\nu} \in \R^{2m}$ corresponds to i.\,i.\,d.\ zero-mean normal noise with variance $\sigma^2 = 10^{-3}$ ($\mathrm{SNR} \approx 60$) and $G = F \otimes F$, where $F$ denotes the one-dimensional discrete Fourier transform with missing frequencies, which we impose to mimic the situation where the system is under-determined and some data must for some reason be discarded. 
The removed frequencies were determined by sampling $100$ logarithmically spaced integers between $10$ and $200$. 
Finally, because the image is piecewise constant, we used first-order TV-regularization. 

\cref{fig:MRI_2d_ML} shows the maximum likelihood (ML) estimate of the image, which is obtained by maximizing the likelihood function $p(\mathbf{x}|\mathbf{y})$. 
In this case, the ML estimate is the same as the least squares (LS) solution of the linear system \cref{eq:data_model_Fourier}. 
\cref{fig:MRI_2d_ADMM,fig:MRI_2d_BCD} illustrate the reconstructions obtained by applying ADMM to the $\ell^1$-regularized inverse problem \cref{eq:l1_RIP} and the SBL-based BCD algorithm. 
The regularization parameter in \cref{eq:l1_RIP} was again fine-tuned by hand and chosen as $\lambda = 4 \sigma^2$. 
While the reconstructions in \cref{fig:MRI_2d_ADMM,fig:MRI_2d_BCD} are comparable, it is important to point out that we did not use any prior knowledge about the noise variance or perform any parameter tuning for the BCD algorithm.

\subsection{Data fusion} 

As a final example we consider a data fusion problem to demonstrate the possible advantage of using the generalized noise model discussed in \cref{sub:likelihood}. 
Recall the piecewise constant signal discussed in \cref{sub:deconvolution_1d}, and assume we want to reconstruct the values of this signal at $n=40$ equidistant grid points, denoted by $\mathbf{x}$.
We are given two sets of data:  
$\mathbf{y}^{(1)}$ corresponds to direct observations taken at $36$ randomly selected locations with added i.\,i.\,d.\ zero-mean normal noise $\boldsymbol{\nu}^{(1)}$ with variance $\sigma_1^2 = 5 \cdot 10^{-1}$, and $\mathbf{y}^{(2)}$ corresponds to blurred observations at $24$ randomly selected locations with added i.\,i.\,d.\ zero-mean normal noise $\boldsymbol{\nu}^{(2)}$ with variance $\sigma_2^2 = 10^{-2}$. 
The blurring is again modeled using \cref{eq:disc_convolution} with a Gaussian convolution kernel and convolution parameter $\gamma = 3 \cdot 10^{-2}$. 
Further, a first-order TV-regularization operator is employed to promote a piecewise constant reconstruction.  

\begin{figure}[tb]
	\centering
  	\begin{subfigure}[b]{0.45\textwidth}
		\includegraphics[width=\textwidth]{%
      		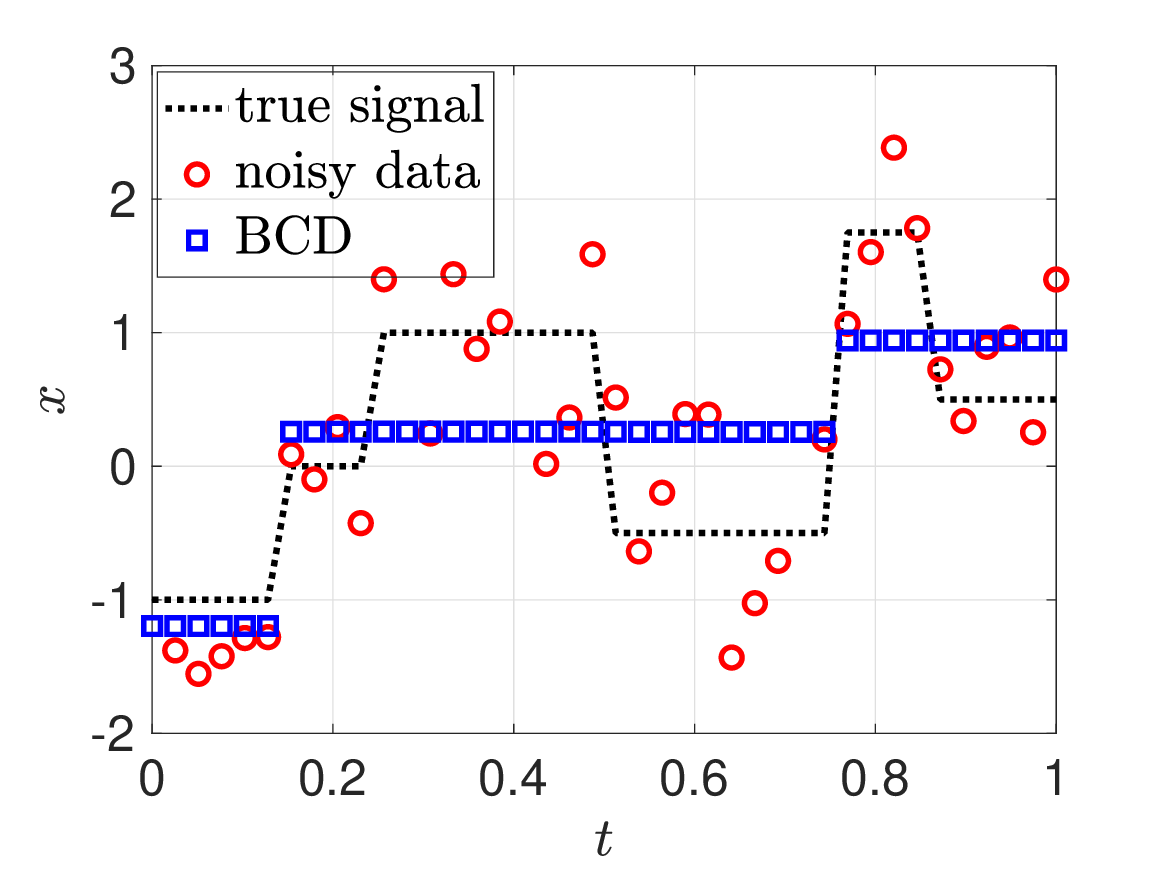} 
    	\caption{Noisy data}
    	\label{fig:dataFusion_1d_sensor1}
  	\end{subfigure}%
  	~
  	\begin{subfigure}[b]{0.45\textwidth}
		\includegraphics[width=\textwidth]{%
      		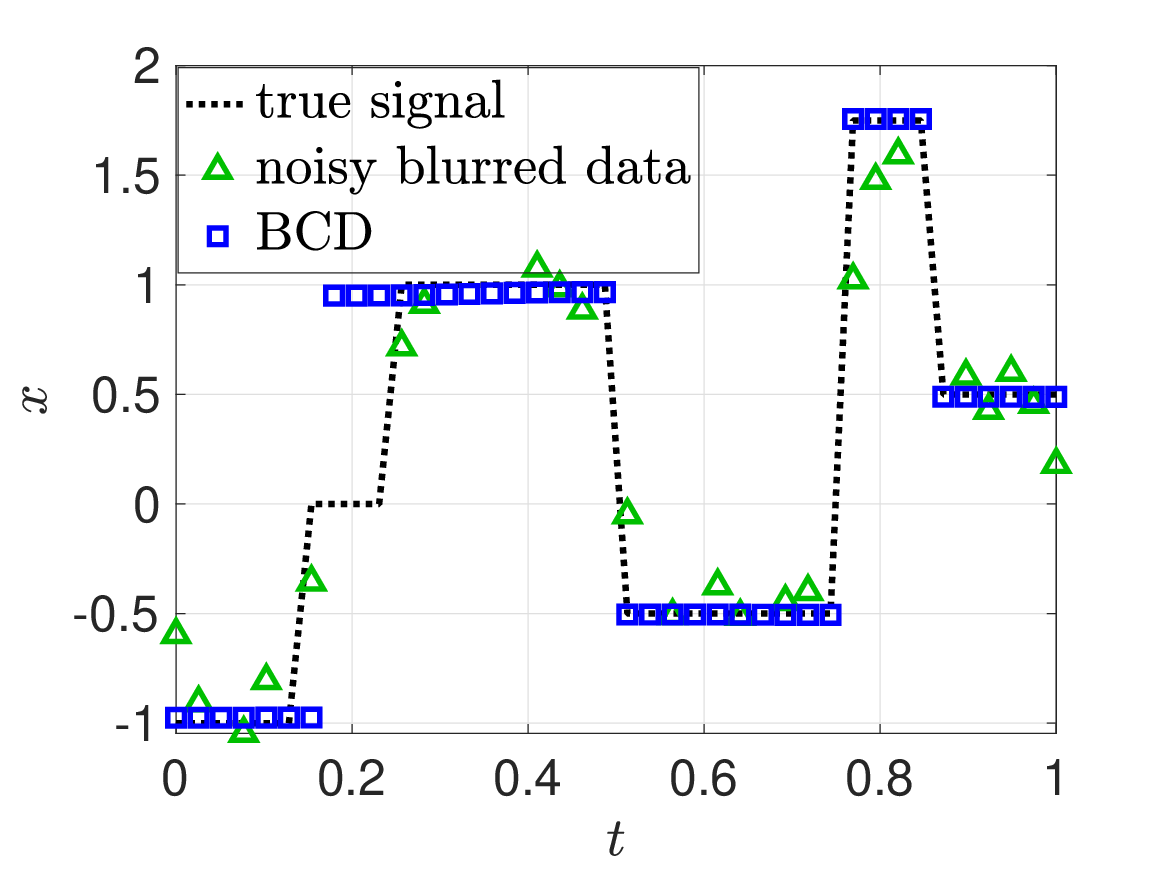} 
    	\caption{Noisy blurred data}
    	\label{fig:dataFusion_1d_sensor2}
  	\end{subfigure}%
	\\ 
	\begin{subfigure}[b]{0.45\textwidth}
		\includegraphics[width=\textwidth]{%
      		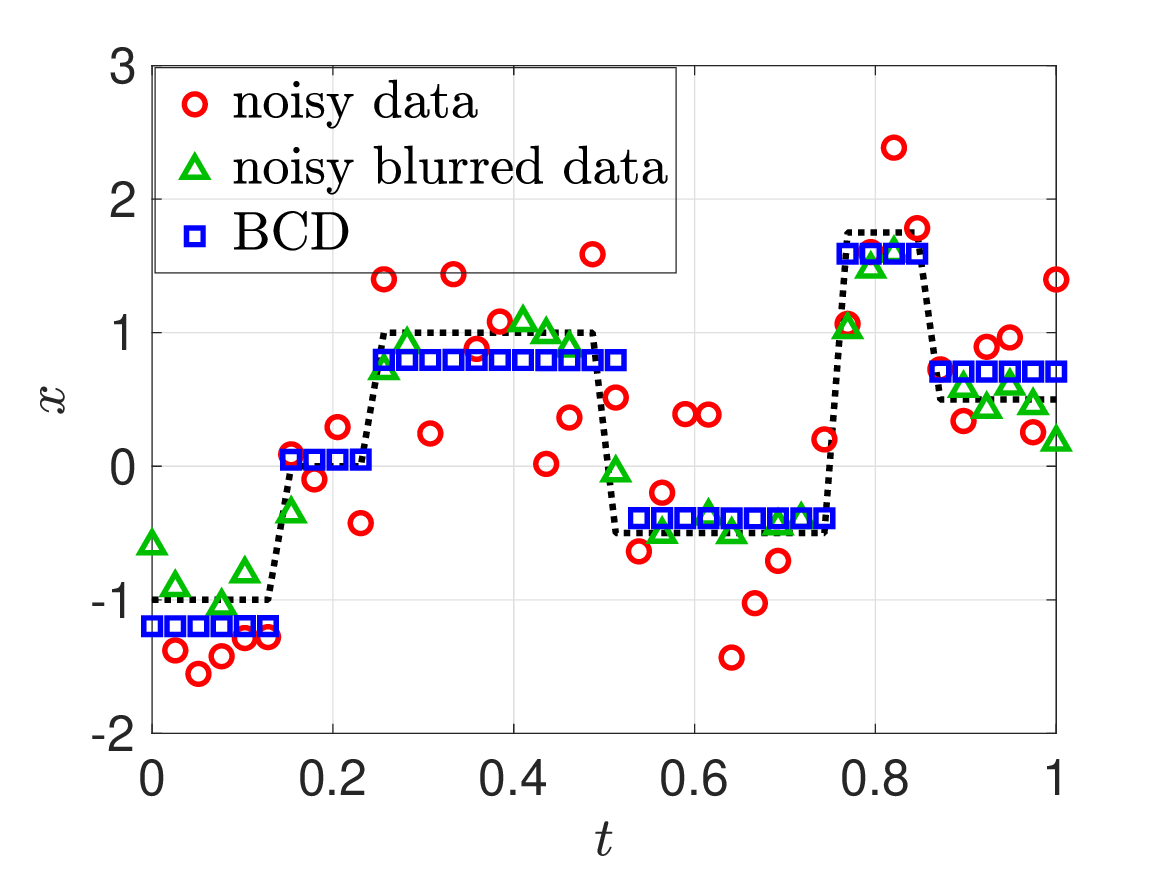} 
    	\caption{Combined, i.\,i.\,d.\ assumption}
    	\label{fig:dataFusion_1d_combined_iid}
  	\end{subfigure}%
  	~
  	\begin{subfigure}[b]{0.45\textwidth}
		\includegraphics[width=\textwidth]{%
      		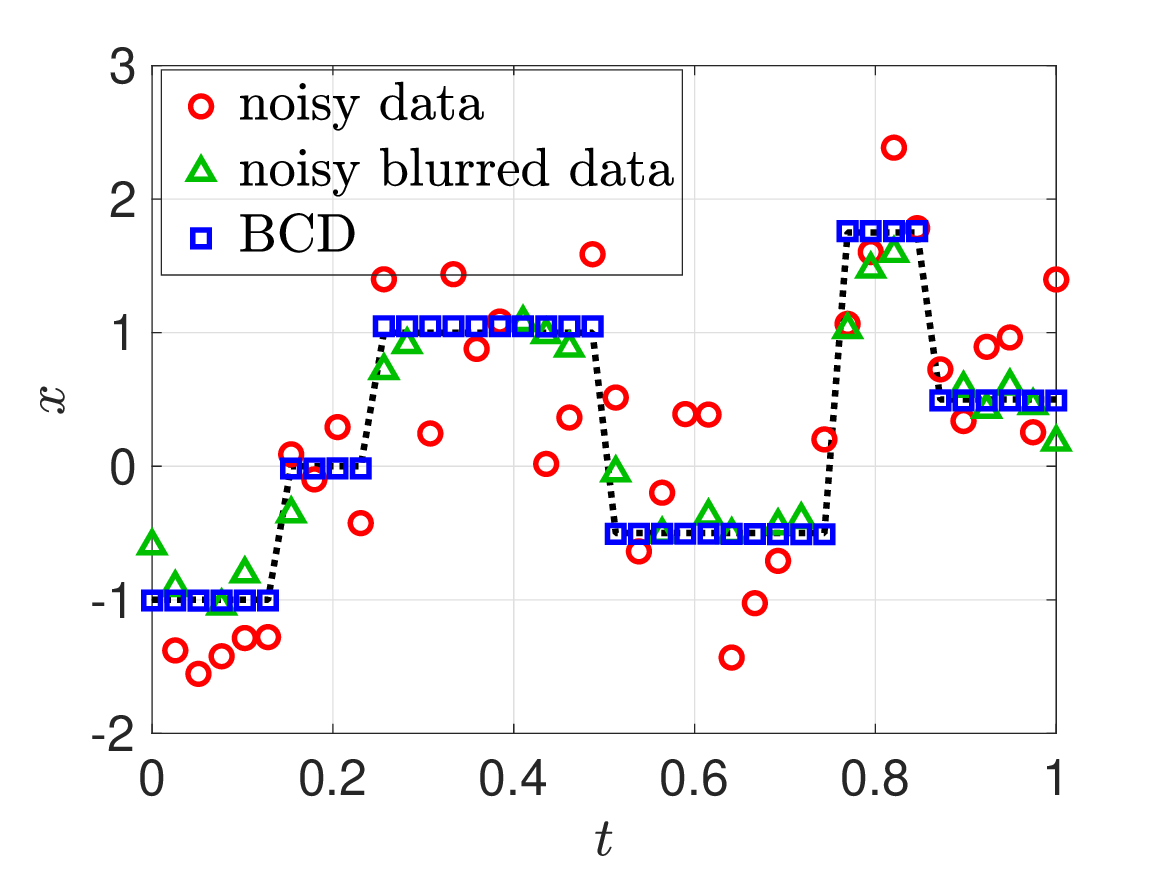} 
    		\caption{Combined, generalized data model}
    		\label{fig:dataFusion_1d_combined}
  	\end{subfigure}%
  	\caption{ 
  	Data fusion example with incomplete noisy and incomplete noisy blurred data. 
  	Top row: Separate reconstructions using the SBL-based BCD algorithm. 
  	Bottom row: Combined reconstructions using the SBL-based BCD algorithm with i.\,i.\,d.\ assumption and using a generalized data model. 
  	}
  	\label{fig:dataFusion_1d}
\end{figure} 

The separate reconstructions by the SBL-based BCD algorithm can be found in \cref{fig:dataFusion_1d_sensor1,fig:dataFusion_1d_sensor2}. 
Both reconstructions are of poor quality, which is due to the high noise variance in the case of $\mathbf{y}^{(1)}$ and to the missing information in the case of $\mathbf{y}^{(2)}$. 
In fact, the reconstruction illustrated in \cref{fig:dataFusion_1d_sensor2} is of reasonable quality except for the region around $t=0.2$, where a void of observations causes the reconstruction to miss the jumps at $t=0.15$ and $t=0.25$. 

Following \cref{expl:data_fusion}, we now fuse the two data sets by considering the joint data model 
\begin{equation}\label{eq:data_model_fusion_tests} 
    \underbrace{%
    \begin{bmatrix} \mathbf{y}^{(1)} \\ \mathbf{y}^{(2)} \end{bmatrix}%
    }_{= \mathbf{y}}
    = 
    \underbrace{%
    \begin{bmatrix} F^{(1)} \\ F^{(2)} \end{bmatrix}%
    }_{= F}
    \mathbf{x}
    + 
    \underbrace{%
    \begin{bmatrix} \boldsymbol{\nu}^{(1)} \\ \boldsymbol{\nu}^{(2)} \end{bmatrix}%
    }_{= \boldsymbol{\nu}},
\end{equation}
where $F^{(1)}$ and $F^{(2)}$ are the forward models describing how $\mathbf{x}$ is mapped to $\mathbf{y}^{(1)}$ and $\mathbf{y}^{(2)}$, respectively. 
Employing the usual likelihood function \cref{eq:likelihood_iid} would correspond to assuming that all the components of stacked noise vector $\boldsymbol{\nu}$ are i.\,i.\,d., which is not true for this example. 
The resulting reconstruction by the BCD algorithm can be found in \cref{fig:dataFusion_1d_combined_iid}. 
In contrast, utilizing the generalized likelihood function \cref{eq:likelihood} with 
\begin{equation}
    A = \diag(\underbrace{\alpha_1,\dots,\alpha_1}_{\text{$m_1$ times}},\underbrace{\alpha_2,\dots,\alpha_2}_{\text{$m_2$ times}}),
\end{equation} 
we can appropriately model that $\boldsymbol{\nu}^{(1)}$ and $\boldsymbol{\nu}^{(2)}$ have different variances. 
The corresponding reconstruction by the BCD algorithm using this generalized data model is provided in \cref{fig:dataFusion_1d_combined}. 
Observe that the reconstruction using the generalized noise model (\cref{fig:dataFusion_1d_combined}) is clearly more accurate than the one for the i.\,i.\,d.\ assumption (\cref{fig:dataFusion_1d_combined_iid}). 
This can be explained by noting that the first data set is larger than the second one, containing $m_1 = 36$ and $m_2 = 24$ observations, respectively. 
At the same time, the data of the first set is less accurate than of the second one, since the variances are $\sigma_1^2 = 5 \cdot 10^{-1}$ and $\sigma_2^2 = 10^{-2}$ ($\mathrm{SNR}_1 \approx 1.6$ and $\mathrm{SNR}_2 \approx 80$), respectively. 
Hence, when using the usual i.\,i.\,d.\ assumption, the first data set $\mathbf{y}^{(1)}$, which is larger but less accurate, more strongly influences the reconstruction than second data set, which is smaller but more accurate. 
Using the generalized data model, on the other hand, the BCD algorithm is able to more appropriately balance the influence of the different data sets. 
 
\section{Concluding remarks} 
\label{sec:summary} 

This paper introduced a generalized approach for SBL and an efficient realization of it by the newly proposed BCD algorithm. 
In contrast to existing SBL methods, we are able to use any regularization operator $R$ as long as the common kernel condition \cref{eq:common_kernel} is satisfied, a usual assumption in regularized inverse problems. 
Further, the BCD algorithm provides us with the full solution posterior $p(\mathbf{x}|\mathbf{y})$ for fixed hyper-parameters rather than just resulting in a point estimate, while being easy to implement and highly efficient.   
Future work will elaborate on sampling based methods for Bayesian inference \cite{bardsley2012mcmc}, which might be computationally more expensive but would also allow sampling from the full joint posterior $p(\mathbf{x},\boldsymbol{\alpha},\boldsymbol{\beta}|\mathbf{y})$. 
This has been addressed to some extent in \cite[Section 6]{calvetti2008hypermodels} for uncertainty quantification in regions of interest. 
Other research directions might include data-informed choices for the parameters $c$ and $d$ in \cref{eq:hyper_priors} and data fusion applications. 
Finally, it would be of interest to combine the proposed generalized SBL framework with generalized Gamma distributions as hyper-priors \cite{calvetti2020sparse} and the hybrid solver from \cite{calvetti2020sparsity}. 

\appendix 

\section{Evidence approach} 
\label{app:improper_evidences} 

In the evidence approach \cite{tipping2001sparse,babacan2010sparse}, the posterior $p( \mathbf{x}, \boldsymbol{\alpha}, \boldsymbol{\beta} | \mathbf{y} )$ is decomposed as 
\begin{equation}\label{eq:BI_evi_decomp}
	p( \mathbf{x}, \boldsymbol{\alpha}, \boldsymbol{\beta} | \mathbf{y} ) 
		= p( \mathbf{x} | \mathbf{y}, \boldsymbol{\alpha}, \boldsymbol{\beta} ) 
			p( \boldsymbol{\alpha}, \boldsymbol{\beta} | \mathbf{y} ).
\end{equation} 
The variables $\mathbf{x}$, $\boldsymbol{\alpha}$, and $\boldsymbol{\beta}$ are then alternatingly updated, with the hyper-parameters $\boldsymbol{\alpha}$ and $\boldsymbol{\beta}$ calculated as the mode (most probable value) of the hyper-parameter posterior $p( \boldsymbol{\alpha}, \boldsymbol{\beta} | \mathbf{y})$. 
By Bayes' law, one has 
\begin{equation}\label{eq:BI_evidence}
	p( \boldsymbol{\alpha}, \boldsymbol{\beta} | \mathbf{y}) 
		\propto  p(\boldsymbol{\alpha}) 
			p(\boldsymbol{\beta}) 
			p( \mathbf{y} | \boldsymbol{\alpha}, \boldsymbol{\beta} ),
\end{equation} 
where the \emph{evidence} $p( \mathbf{y} | \boldsymbol{\alpha}, \boldsymbol{\beta} )$ can be determined by marginalizing out the unknown solution $\mathbf{x}$, which yields  
\begin{equation}\label{eq:evidence_marg}
	p( \mathbf{y} | \boldsymbol{\alpha}, \boldsymbol{\beta} ) 
		= \int p( \mathbf{y} | \mathbf{x}, \boldsymbol{\alpha} ) p( \mathbf{x} | \boldsymbol{\beta} ) \intd \mathbf{x}. 
\end{equation} 
Some basic but lengthy computations are then used to obtain  
\begin{equation} 
\begin{aligned}
    p( \mathbf{y} | \boldsymbol{\alpha}, \boldsymbol{\beta} ) 
        = (2 \pi)^{(n-m)/2} \det(A)^{1/2} \det(B)^{1/2} \det(C)^{-1/2}  
        \exp \left\{ -\frac{1}{2} \mathbf{y} \Sigma^{-1} \mathbf{y} \right\},
\end{aligned}
\end{equation} 
where $\Sigma = A^{-1} + F ( R^T B R )^{-1} F^T$.
Also see \cite[Section 3]{babacan2010sparse}.
However, this assumes that $R^T B R$ is invertible, which is not the case whenever $\kernel(R) \neq \{\mathbf{0}\}$. 

\section*{Acknowledgements}
This work was supported by AFOSR \#F9550-18-1-0316 (Glaubitz and Gelb), NSF-DMS \#1502640, NSF-DMS \#1912685, ONR \#N00014-20-1-2595 (Gelb), NSF-DMS \#1521661, and NSF-DMS \#1939203 (Song).

\bibliographystyle{siamplain}
\bibliography{references}

\end{document}